\documentclass[twocolumn]{autart}
\usepackage{generic}
\usepackage{cite}
\usepackage{amsmath,amssymb,amsfonts}
\usepackage{algorithmic}
\usepackage{graphicx}
\usepackage{algorithm,algorithmic}
\usepackage{hyperref}
\hypersetup{hidelinks=true}
\usepackage{textcomp}

\usepackage{cite}
\usepackage{amsmath,amssymb,amsfonts}
\usepackage{dsfont}
\usepackage{algorithmic}
\usepackage{graphicx}
\usepackage{textcomp}
\usepackage{xcolor}

 \usepackage{amsmath}
\usepackage{amsfonts}
\usepackage{amssymb}
\usepackage{mathtools}
\usepackage{graphicx}
\usepackage{hyperref}
\usepackage{setspace}

\newtheorem{theorem}{Theorem}[section]
\newtheorem{corollary}{Corollary}[theorem]
\newtheorem{lemma}[theorem]{Lemma}

\newtheorem{definition}{Definition}[section]

\newcommand{\bE}{\mathbb{E}}
\newcommand{\bsx}{\boldsymbol{x}}

\newcommand{\tbsx}{\tilde{\boldsymbol{x}}}
\newcommand{\bsu}{\boldsymbol{u}}
\newcommand{\bsy}{\boldsymbol{y}}
\newcommand{\bsv}{\boldsymbol{v}}
\newcommand{\bsw}{\boldsymbol{w}}
\newcommand{\bW}{\boldsymbol{W}}
\newcommand{\bbW}{\mathbb{W}}
\newcommand{\bM}{\boldsymbol{M}}
\newcommand{\bbM}{\mathbb{M}}

\newcommand{\bA}{\boldsymbol{A}}
\newcommand{\bbA}{\mathbb{A}}

\newcommand{\bB}{\boldsymbol{B}}
\newcommand{\bbB}{\mathbb{B}}
\newcommand{\bQ}{\mathbf{Q}}
\newcommand{\bR}{\mathbf{R}}
\newcommand{\bbR}{\mathbb{R}}
\newcommand{\calR}{\mathcal{R}}

\newcommand{\bbS}{\mathbb{S}}

\newcommand{\bbK}{\mathbb{K}}

\newcommand{\bbI}{\mathbb{I}}

\newcommand{\op}{\text{op}}

\newcommand{\bsxf}{\boldsymbol{x}^{\boldsymbol{f}}}
\newcommand{\bsxp}{\Breve{\boldsymbol{x}}}

\newcommand{\bsuf}{\boldsymbol{u}^{\boldsymbol{f}}}
\newcommand{\bsup}{\Breve{\boldsymbol{u}}}

\newcommand{\bswf}{\boldsymbol{w}^{\boldsymbol{f}}}
\newcommand{\bswp}{\Breve{\boldsymbol{w}}}

\newcommand{\calS}{\mathcal{S}}
\newcommand{\bsf}{{\boldsymbol{f}}}

\begin{document}
\begin{frontmatter}

\title{Quadratic Optimal Control of Graphon Q-noise Linear Systems }
\thanks{The authors are with the Department of Electrical and Computer Engineering, McGill University, Montreal, QC H3A 0E9, Canada (e-mail:alex.dunyak@mail.mcgill.ca; peterc@cim.mcgill.ca).}
	\thanks{Work supported by   *A. Dunyak  NSERC 2019-05336, **P.E. Caines:  NSERC 2019-05336, ARL W911NF1910110, Air Force OSR Grant FA9550-23-1-0015 }

\author{Alex Dunyak and Peter E. Caines}
    \begin{abstract}
       The modelling of linear quadratic Gaussian optimal control problems on large complex networks is intractable computationally. Graphon theory provides an approach to overcome these issues by defining limit objects for infinite sequences of graphs permitting one to approximate arbitrarily large networks by infinite dimensional operators. This is extended to stochastic systems by the use of Q-noise, a generalization of Wiener processes in finite dimensional spaces to processes in function spaces. The optimal control of linear quadratic problems on graphon systems with Q-noise disturbances are defined and shown to be the limit of the corresponding finite graph optimal control problem. The theory is extended to low rank systems, and a fully worked special case is presented. In addition, the worst-case long-range average and infinite horizon discounted optimal control performance with respect to Q-noise distribution are computed for a small set of standard graphon limits.
    \end{abstract}
\end{frontmatter}
\section{Introduction}
    Large graphs are common objects in contemporary systems modelling and analysis, in particular for the purposes of optimization and control. Indeed, from the internet to  electrical generation and distribution to social networks, complex networks have been a focus of research for decades. However, global modelling and analysis  problems  are intractable with standard methods for all sufficiently large networks.   
        
    One approach to handling large networks is to use the theory of graphons \cite{lovasz_large_2012}.  Informally, a graphon is a function on the unit square which represents a limit of the adjacency matrices of a sequence of graphs. Consequently, using graphons  for modelling large systems allows for the approximation  of very large network within a functional analysis framework and hence enables their modelling, analysis and design.
   
    The use of graphons in system dynamics was initiated by Medvedev \cite{medvedev_nonlinear_2014}. Previous work on control via graphons has been primarily concerned with deterministic  systems (\cite{gao_spectral_2019}\cite{gao_graphon_2020}), while stochastic mean field games on graphons have been investigated in \cite{caines_graphon_2021}, \cite{caines_graphon_2018}, \cite{tchuendom_master_2020}, \cite{gao_lqg_2021}, \cite{lacker_label-state_2023} and with Q-Wiener processes in \cite{liu_hilbert_2024}. Graphon Mean Field Games have been extended in particular to the case of large sparse graphs in \cite{lacker_case_2022}, \cite{cui_learning_2022}, and \cite{caines_embedded_2022}. The use of Q-Wiener processes is an attempt to solve the measurability problem of applying idiosyncratic Wiener processes at each node, which was addressed by Aurell et al. for linear quadratic graphon mean field games using the Fubini extension in \cite{aurell_stochastic_2022} and for epidemic games in \cite{aurell_finite_2022}.
    
    In Dunyak and Caines \cite{dunyak_linear_2022}, space-time Gaussian noise on the unit interval (\cite{gawarecki_stochastic_2011} \cite{fabbri_stochastic_2017} \cite{kukush_gaussian_2019}), termed Q-noise, was introduced as a limit object for sequences of systems on graphs with Brownian disturbances. Medvedev and Simpson \cite{medvedev_numerical_2023} presented a numerical method of simulating such systems. This article demonstrates that linear quadratic Q-noise optimal control problems on large graphs can be approximately solved by the graphon limit of their system. In this article, Linear Quadratic Gaussian problem in Hilbert spaces are approximately solved are approximately solved by analysis of the graphon limits.  The analysis relies on the Hilbert space methods in Ichikawa \cite{ichikawa_dynamic_1979} which is extended here to long-range average and infinite discounted cost LQG problems, the former of which is derived with a limit argument and the latter of which is solved via the corresponding algebraic Riccati equation in the case of purely local controls.

    In the following subsections, we provide the motivation for modelling linear quadratic Gaussian systems on large networks with their graphon and Q-noise limits. In Section II, we define the notation used in this article, as well as summarize relevant prior results for Q-noise systems. In Section III, we present the formal proof that the finite dimensional linear quadratic Gaussian system converges to the infinite dimensional linear quadratic Q-noise system, as well as presenting the long-range average and exponential discounting problems. Section IV extends the analysis of low rank graphons presented in  \cite{gao_subspace_2021} to Q-noise systems. Section V demonstrates the utility of the Q-noise approach in that the solution to an LQG system problems on large unweighted random graphs are shown to be well-approximated by lower-order systems derived from the graphon limit of the original system. Section VI provides some directions for future research.
    
    \subsection{Motivation: Networked Systems and Graphons}
    Define two graphs $G^N_A = (V_N,E^N_A)$ and $G^N_B = (V_N,E^N_B)$ with $N < \infty$ vertices, with associated adjacency matrices $A^N$ and $B^N$. Let $x^N:[0,T]^N \to \mathbb{R}$ be a vector of states, where the $i$th value is associated with the $i$th vertex of the graph, and let $u^N:[0,T]^N \to \mathbb{R}^N$ be the control input at each vertex. For clarity of notation, systems where each node has a scalar state are considered below. The theory extends to systems with vector states at each node in a straightforward manner. Let the $(i,j)$th entry of the matrices $A^N$ and $B^N$ represent the impact of the state and control input at node $i$ on node $j$, respectively. For each node, define a Brownian motion such that the $N$-component Brownian motion $W^N$ has strictly positive covariance matrix $Q^N$. Let $a_N$ and $b_N$ be constants describing the impact of the state of a node and its control on itself.
	
	Finally define a network-averaged control system \cite{gao_graphon_2020} on a graph with the following equation for each node,
	\begin{align}
	    dx^i_t =& (a_N x^i_t +\sum_{j=1}^N A^N_{ij} x_t^j \nonumber\\
	     &~~+ b_N u^i_t +\sum_{j=1}^N B^N_{ij} u_t^j)dt + dW^N_i(t),
	\end{align}
	or in vector form,
	\begin{align} \label{graph_system}
	    dx_t^N = ((A^N + a_N) x^N_t + (B^N + b_N) u^N_t)dt + dW^N(t).
	\end{align}
    Subject only to the assumption that the entries of $A^N$ and $B^N$ are uniformly bounded in $N$, the sequences of adjacency matrices $\{A^N\}$ and $\{B^N\}$, $1 \leq N < \infty$, defined on the unit square converge, as $N$ tends to infinity, to their (not necessarily unique) associated graphon limits \cite{lovasz_large_2012}, which are bounded measurable functions mapping $[0,1] \times [0,1] \to [0,1]$. These are denoted $\bA$ and $\bB$ (as in \cite{gao_graphon_2020})\footnote{In order to disambiguate the convergence of the adjacency matrix $A^N$ to the graphon $\bA$, the scaling term $\frac{1}{N}$ is omitted when $A^N$ is acting as an operator. This scaling is to ensure that the summation $A^N x^N$ is bounded and converges to the correct integral.}. When the underlying graph is undirected, its graphon is also symmetric. Denote the graphon limit system as
    \begin{align}
        d\bsx_t = ( (\bA + a\bbI)\bsx_t + (\bB + b \bbI) \bsu_t ) + d\bsw_t,
    \end{align}
    where $\bsx_t$ and $\bsu_t$ are square-integrable functions on the unit interval, $\bA$ and $\bB$ are graphons, $a$ and $b$ are real constants, $\bbI$ is the identity operator, and $\bsw_t$ is a Q-noise, a generalization of Gaussian noise from finite-dimensional vectors to random time-varying functions the unit interval as defined in Section \ref{sec:Q-noise}, together with the conditions in Theorem \ref{thm:mean_convergence} for the existence of the limit in mean-square.
    
\subsection{Linear Quadratic Q-noise Control}

    Let $x_t^N$ be the state of a networked control system on a graph $G^N$ as given in equation (\ref{graph_system}). Suppose that $M^N$ is an $N \times N$ positive matrix, and $R^N$ is an $N\times N$ strictly positive matrix. Then the associated linear quadratic Gaussian optimal control problem for a control system with terminal time $T$ is defined via the infimization of the performance function:
    \begin{align}
        \inf_{u_0^T} J(x_0,u) = \inf_{u_0^T} \bE\left[\int_0^T x^{N*}_t M^N x^N_t + u_t^{N*} R^N u_t^N dt\right].
    \end{align}
    
    The solution of the limit problem takes the same form as the standard finite dimensional LQG problem, but the equations have operator valued coefficients and the solutions are operator valued.  This work analyses the properties of the operator  limits of such sequences of network averaged optimal control problems and it is shown that the solutions of the limit problems are obtained via the operator limits of the associated Riccati equations.

\subsection{The Special Case of Finite Rank Systems}

    The systems described in the prior section are defined in the space of square-integrable functions on the unit interval, $L_2[0,1]$. In general, it is not possible to find a closed form solution for such a system. However, when a system's associated graphon parameters and Q-noise covariance function are finite-rank, then the state of the system evolves on a finite dimensional space. This is explored in Section \ref{sec:low_rank}.



\section{Preliminaries}
	
	\subsection{Notation}

\begin{itemize}
        \item The set of vectors of real numbers of dimension $m$ is denoted $\calR^m$.
        \item Graphons (i.e. bounded symmetric $[0,1]^2$ functions used as the kernels of linear integral operators) are denoted in italicized bold capital letters, such as $\bA$, $\bB$, and $\bM$. 
        \item $L_2[0,1]$ denotes the Hilbert space of real square-integrable functions on the unit interval.
        In addition, $L_2[0,1]$ is equipped with the standard inner product, denoted $\langle \bsu, \bsv \rangle$.
        For any function $\bsv$, $\bsv^*$ denotes the adjoint of $\bsv$. As such, $\langle \bsu, \bsv \rangle$ is sometimes written as $\bsv^* \bsu$.
        \item The identity operator in both $L_2[0,1]$ and finite dimensional spaces is denoted $\bbI$.
        \item Operators of the form $\bbA$ and $\bbB$ have the structure $\bbA = \bA + a \bbI$, where $\bA$ is a graphon and $a$ is a real scalar. Let $\mathcal{M}$ denote the set of these operators.
        \item A linear integral operator with the kernel $\bQ: [0,1]^2 \to \mathcal{R}$ acting on a function $\bsf \in L_2[0,1]$ is defined by
    	\begin{align}
    		(\bQ  \bsf)(x) = \int_0^1 \bQ(x,y) \bsf(y) dy, \quad \forall ~ x \in [0,1].
    	\end{align}
        \item The operators $\bQ$ are equipped with the standard operator norm $||\bQ||_\op$. When unambiguous, the argument is dropped. 
        \item A symmetric function $\bQ : [0,1]^2 \to \mathcal{R}$ is non-negative if the following inequality is satisfied for every function $\bsf \in L_2[0,1]$,
    	\begin{align} \label{notation:positive}
        	0 &\leq 	\int_0^1 \int_0^1 \bQ (x,y)\bsf^*(x)\bsf(y) dxdy  \\
            & := \langle \bQ \bsf, \bsf \rangle <\infty. \nonumber
    	\end{align}
        Additionally, denote $\mathcal{Q}$ to be the set of bounded symmetric non-negative functions. The set $\mathcal{Q}$ serves as the set of valid covariance functions for the class of stochastic processes called Q-noise processes.
        \item For a given linear system $\bsx_t$ satisfying $\dot{\bsx}_t = \bA_t \bsx_t$, $\Phi(t,s)$ denotes the semigroup solution operator solving $\bsx_t = \Phi(t,s)\bsx_s$ for any given initial condition $\bsx_s$. 
        \item Q-noise processes (stochastic processes over the time interval $[0,T]$) will be denoted as $\bsw_t$. For each $t \in [0,T]$, $\bsw_t$ is an $L_2[0,1]$ function. The precise definition of a Q-noise process is given in Section \ref{sec:Q-noise}.
        \item A partition of the unit interval of $N$ increments is denoted $P^N = \{P_1, \cdots, P_N\}$, where $P_1 = [0,\frac{1}{N}]$ and $P_i = (\frac{i-1}{N},\frac{i}{N}]$. An $L_2[0,1]$ function which is piece-wise constant on the unit interval is denoted $\bsv^{[N]}$, and a self-adjoint $L_2[0,1]$ operator $\bM$ which is piece-wise constant on the Cartesian product $P^N \times P^N$ is denoted $\bM^{[N]}$ (or $\bbM^{[N]}$, if it is of the form $\bbM^{[N]} = \bM^{[N]} + m\bbI$). This formulation is necessary for mapping $N \times N$ adjacency matrices of networks to functions on the unit square, as in Section \ref{sec:networks}.
    \end{itemize}



\subsection{Q-noise Axioms} \label{sec:Q-noise}
Q-noise processes, first applied to graphon systems in \cite{dunyak_linear_2022}, are $L_2[0,1]$ valued random processes that satisfy the following axioms.
\begin{enumerate}
		\item Let $\bQ \in \mathcal{Q}$, and let ($[0,1] \times [0,T] \times \Omega, \mathcal{B}([0,1] \times [0,T] \times \Omega),\mathbb{P})$ be a probability space with the measurable random variable $\bsw (\alpha, t,\omega):[0,1]  \times [0,T] \times \Omega\to \mathcal{R}$ for all $t \in [0,T], \alpha \in [0,1]$, $\omega \in \Omega$. For notation, $\omega$ is suppressed when the meaning is clear.
		\item For all  $\alpha \in [0,1]$, $\bsw(\alpha,t) - \bsw(\alpha,s)$ is a Wiener process increment in time for all $t,s \in [0,T]$, with $\bsw(\alpha,t) - \bsw(\alpha,s) \sim \mathcal{N}(0, |t-s| \bQ(\alpha,\alpha))$ where $\bsw(\alpha,0) = 0$ for all $\alpha\in [0,1]$.
		\item Let $\bsw_{t-t'}(\alpha) = \bsw(\alpha,t) - \bsw(\alpha,t')$. Then, 
        \begin{align*}
            \mathbb{E}[\bsw_{t-t'}(\alpha) \bsw_{s-s'}(\beta)] =|[t , t'] \cap [s,s']|\cdot \bQ(\alpha,\beta).
        \end{align*}
		\item For almost all $s,t \in [0,T]$, $\alpha, \beta \in [0,1]$, and $\omega \in \Omega$, $\bsw(\alpha,t,\omega) - \bsw(\beta,s,\omega)$ is Bochner-integrable\footnote{The Bochner integral of a random variable $X:(\Omega,\mathcal{B},\mu) \to E$ where $E$ is a Banach space is defined as the limit of the sum of simple functions taking a finite set of values $X_n(\omega)$, analogously to Lebesgue integration \cite[Sec. 1.1.3]{fabbri_stochastic_2017}. } as a function taking values in the Banach space of a.s. piece-wise continuous functions of $[s,t] \in [0,T]$.
	\end{enumerate}

\textit{An orthonormal basis example:}
    Let $\{W_1,W_2,\cdots\}$ be a sequence of independent Brownian motions. Let $\bQ \in \mathcal{Q}$ have a diagonalizing orthonormal basis $\{\phi_k\}_{k=1}^\infty$ with eigenvalues $\{\lambda_k\}_{k=1}^\infty$. Then
    \begin{equation}
        g(\alpha,t,\omega) = \sum_{k=1}^\infty \sqrt{\lambda_k}\phi_k(\alpha)W_k(t,\omega)
    \end{equation}
    is a Q-noise process. The common name for this formulation in the literature is $Q$-Wiener process (\cite{gawarecki_stochastic_2011} \cite{fabbri_stochastic_2017}).
    



\subsection{Operators on Q-noise Noise}

\begin{definition}
      $\mathcal{M}$ shall denote the set of operators of the form $\bbM = \bW + c\bbI$, where $\bW$ is a bounded self-adjoint Hilbert-Schmidt integral operator (hence possessing square-summable eigenvalues) mapping $\mathcal{L}^2[0,1]$ to $\mathcal{L}^2[0,1]$, $c>0$ is a positive constant, and $\bbI$ is the identity operator on $\mathcal{L}^2[0,1]$.
	\end{definition}
    \begin{definition}[Operators on $Q$-Space Noise]
		Let $\bQ \in \mathcal{Q}$ and $\bsw_t$ be a $Q$-space noise. Let $\bbM \in \mathcal{M}$, and let $s < t \in [0,T]$. Then the action of $\bbM$ on $\bsw_{t-s}(\cdot) :=  \bsw(\cdot, t) - \bsw(\cdot, s)$ is defined by the following Lebesgue integral for $s,t \in[0,T], ~ \alpha \in [0,1]$,
		\begin{align}
			(\bbM & \bsw_{t-s}(\cdot)) (\alpha) = \int_0^1 \bbM(\alpha,z) \bsw_{t-s}(z)dz.
		\end{align}
	\end{definition}
    \begin{lemma}\label{Operator:centered}
	Let $\bbM = \bW + c\bbI\in \mathcal{M}$. Then	$(\bbM \bsw_{t-s}) (\alpha)$ is a centered random variable for all $\alpha \in [0,1]$ and $s,t \in [0,T]$.
	\end{lemma}
    \textit{Proof:} As $\bW$ is a bounded operator, $\bE[\bsw_{t-s}] = 0$ by Axiom 2, and since $\bsw_{t-s}$ is assumed to be Bochner-integrable by Axiom 4, the expected value is 
    \begin{align}
        &\bE[\bbM  \bsw_{t-s}](\alpha) = \bE[\bW \bsw_{t-s}](\alpha) + c \bE[\bsw_{t-s}](\alpha) \nonumber \\
        &= (\bW + c\bbI) \bE[ \bsw_{t-s}](\alpha)  = 0.  \quad \square
    \end{align}
	
	By associativity, an operator $\bbW$ acting on an operator $\bbM$ acting on a function $\bsx$ is denoted $(\bbW (\bbM \bsx)) (\alpha) = (\bbW \bbM  \bsx) (\alpha)$ when the following iterated integral exists,
	\begin{equation}
		((\bbW  \bbM)  \bsx )(\alpha) = \int_{0}^{1} \bbW(\alpha,z) \int_0^1 \bbM(z,y) \bsx(y) dy dz.
	\end{equation}
	\begin{theorem}\label{Q_cov}
		Let $\bsw_{t-t'} (\cdot) = \bsw(\cdot,t) - \bsw(\cdot,t')$ and $\bsw_{s-s'} (\cdot) = \bsw(\cdot,s) - \bsw(\cdot,s')$ be two time increments of a $Q$-space process, $\bQ \in \mathcal{Q}$, and $\bbM \in \mathcal{M}$ with $\bbM \bQ \bbM^* \in \mathcal{Q}$. Then,
		\begin{align}
			\text{cov}((\bbM&  \bsw_{t-t'})(\alpha), (\bbM  \bsw_{s-s'})(\beta)) = \nonumber\\
			& \quad \quad= |[t, t']\cap[s,s']| (\bbM  \bQ  \bbM^*)(\alpha,\beta)
		\end{align}
	\end{theorem}
	\textit{Proof}: Recalling Lemma \ref{Operator:centered} the covariance is given by 
	\begin{align}
		\text{cov}((\bbM&  \bsw_{t-t'})(\alpha), (\bbM  \bsw_{s-s'})(\beta)) = \nonumber\\
		& ~~  =\bE[ (\bbM   \bsw_{t-t'})(\alpha)  ((\bsw_{s-s'})^*  \bbM^*)(\beta)]
	\end{align}
	As $\bsw_{t-s}$ is Bochner-integrable by Axiom 4 and $\bbM$ is a bounded self-adjoint operator, the operator $\bbM$ can be exchanged with the expectation,
	\begin{align}
	    \text{cov}&((\bbM  \bsw_{t-t'})(\alpha), (\bbM  \bsw_{s-s'})(\beta)) = \nonumber\\
	    &  = (\bbM  \bE[\bsw_{t-t'}\bsw^*_{s-s'}] \bbM^*)(\alpha,\beta)\\
	    &  = |[t, t']\cap[s,s']|\nonumber\\
	    & ~~~~\cdot\int_0^1 \int_0^1 \bbM(\alpha,y) \bQ(y,z) \bbM^*(z,\beta)dydz. \quad \square \nonumber
	\end{align}



\subsection{Linear Dynamical Control Systems}
    \begin{definition}[Q-noise Dynamical Systems]
		Let $\bsx:[0,1] \times [0,T] \to \mathcal{R}$ be an $\mathcal{L}^2[0,1] \times [0, T]$ function with a given initial condition $\bsx(\cdot,0) = \bsx_0$. Let $\bbA_t, \bbB_t \in \mathcal{M}$ be bounded linear operators from $\mathcal{L}^2[0,1]$ to $\mathcal{L}^2[0,1]$ such that $\bbA_t \bQ \bbA_t^* \in \mathcal{Q}$. This defines a Q-noise denoted $\bsw_t$. Let a control input $\bsu_t:[0,T] \to \mathcal{L}^2[0,1]$ be a function adapted to the filtration $\mathcal{F}_t$, consisting of all measurable functions of the state of the system $\bsx_s, 0 \leq s \leq t$.
		
		Then, a linear dynamical system with Q-noise is an infinite dimensional differential system satisfying the following equation,
		\begin{equation} \label{linear_sys} 
			d \bsx_t(\alpha) = ((\bbA_t \bsx_t)(\alpha) + (\bbB_t \bsu_t)(\alpha))dt + d\bsw(\alpha, t),
		\end{equation}
        where for a partition of $[0,t]$, $[0,t_2, \cdots, t_{N-2},t]$,
		\begin{equation}
			\int_{0}^{t}d\bsw(\alpha,s) = \lim_{N \to \infty} \sum_{k = 1}^N (\bsw(\alpha,t_{k+1}) - \bsw(\alpha,t_k))
		\end{equation}
		in the mean-squared convergence sense.
	\end{definition}
    \begin{definition}[Mild solution]
        A mild solution (see \cite[Sec. 3.1]{gawarecki_stochastic_2011}) to a system $\bsx_t: [0,T] \times \Omega \to L_2[0,1]$ satisfying
        \begin{align}
            d\bsx_t = (\bbA_t \bsx_t + \bbB_t \bsu_t)dt + d\bsw_t, ~ \bsx_0 \in L_2[0,1]
        \end{align} on $[0,T]$ is given by
        \begin{equation}
            \bsx_t = \Phi(t,0)\bsx_0 + \int_0^t \Phi(t,s) \bbB_t \bsu_t ds + \int_0^t \Phi(t,s) d\bsw_s.
        \end{equation}
    \end{definition}
    In the development of this paper, $\bbA$ and $\bbB$ will be taken to be bounded, time-invariant operators. When $\bbA_t$ is constant, then $\Phi(t,s) = e^{\bbA (t-s)}$.
\subsubsection{Motivation for Q-noise models}

To demonstrate why the Q-noise framework is necessary for the modeling of very large linear systems, consider two linear stochastic systems of dimension $N = 300$. For each index $i \in \{1,...,N\}$, the $i$th state of both systems satisfies the following stochastic differential equation with Brownian disturbance $W^i_t$,
\begin{align} \label{eq:indep_noise}
    dx^i_t =& \frac{1}{N}\sum_{j=1}^N \cos\left(\pi \left(\frac{i}{N} - \frac{j}{N}\right)\right) x^j_t dt + dW^i_t, \nonumber\\
    x^i_0 =& 1. 
\end{align}

\begin{figure}
    \centering
    \includegraphics[width=\linewidth]{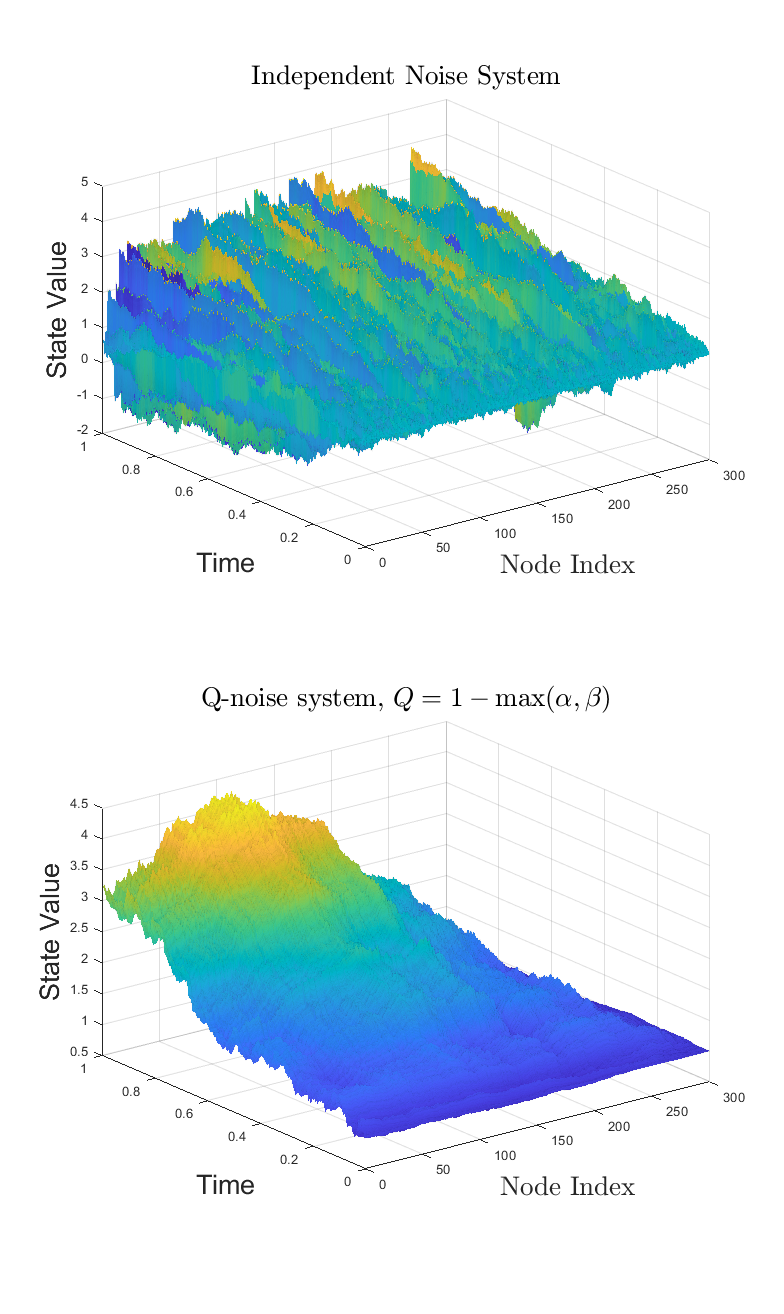}
    \caption{Top: Trajectory of system (\ref{eq:indep_noise}) with independent noise at each node. Because the magnitude of the independent noise is so high, there is no clear system structure in the state trajectory. Bottom: the state trajectory of system (\ref{eq:indep_noise}) with a Q-noise disturbance. While there is clearly noise present in the system, there is an overall structure suggesting that the limit will be continuous in both space and time.}
    \label{fig:indep_noise}
\end{figure}

In the first system, assume that the covariance matrix $Q$ of the Brownian motion disturbance between state $i$ and state $j$ is given entry-wise by $Q_{ij} = 1 - \max\left(\frac{i}{N}, \frac{j}{N}\right)$. This system converges to a Q-noise system in the $L_2[0,1]$ norm as the number of nodes increases. In the second system, assume that the Brownian motion disturbances $W^i_t$ are independent. The result is a system where the noise of individual states overpowers the trajectory of the overall system. The sample trajectories of two such systems with terminal time $T=1$ can be seen in Fig. \ref{fig:indep_noise}.

In addition, there is the fundamental measurability consideration that for real $\alpha \in [0,1]$ parameterized Wiener processes $B_\alpha(t)$ which are independent with respect to $\alpha \in [0,1]$ do not exist as well-defined stochastic processes.


\subsection{Networks} \label{sec:networks}
Consider a network-averaged control system of the form (\ref{graph_system}),
    \begin{align}
	    dx_t = ((A^N x_t + B^N u_t) + \alpha_N x_t + \beta_N u_t)dt + dW^N(t).
	\end{align}
	The finite dimensional system is mapped to piecewise constant functions on the unit square (see \cite{gao_graphon_2020}). Define the uniform partition on the unit interval as $P^N = \{P_1, \cdots, P_N\}$, where $P_1 = [0,\frac{1}{N}]$ and $P_i = (\frac{i-1}{N},\frac{i}{N}]$. Then, the following step function graphon for $N$ nodes corresponding to $A^N$ is defined for all $\alpha,\beta \in [0,1]$:
	
	\begin{equation}
	    \bA^{[N]}(\alpha,\beta) = \sum_{i=1}^N \sum_{j=1}^N A^N_{ij}\mathds{1}_{P_i}(\alpha)\mathds{1}_{P_j}(\beta).
	\end{equation}
	
    A similar function can be defined for $\bB^{[N]}(\cdot,\cdot)$. One can define the step function control $\bsu_t^{[N]}:[0,1] \times [0,T] \to \mathcal{R}$ for such a graphon system using the control $u_t$,
    
    \begin{equation}
        \bsu_t^{[N]}(\alpha) = \sum_{i = 1 }^N \mathds{1}_{P_i}(\alpha)u_t(i) ~~ \alpha \in [0,1].
    \end{equation}
    
    Finally, define the covariance of the disturbance as a piecewise constant function analogously to the finite dimensional adjacency matrix. Let $\bsw^{[N]}:[0,1] \times [0,1] \to \mathcal{R}$ be a Q-noise with covariance defined by the following equation,

    \begin{equation}
        \bQ^{[N]} (\alpha,\beta) = \sum_{i=1}^N \sum_{j=1}^N Q^N_{ij} \mathds{1}_{P_i}(\alpha)\mathds{1}_{P_j}(\beta).
    \end{equation}

    A piecewise constant partition section of the unit interval $\{\bbS^N_k\}_{k=1}^N$ mapping $[0,1] \to \mathcal{R}$ is defined via
    \begin{equation}
        \bbS^N_k(\alpha) := \mathds{1}_{P_k}(\alpha), ~ \alpha \in [0,1].
    \end{equation}
    This defines an orthogonal set corresponding to the standard $\mathcal{R}^N$ basis  $\{e_1, e_2, ..., e_N\}$.
    
	Then the corresponding systems in $\mathcal{L}^2[0,1]$ can be expressed as 
	\begin{align} \label{eq:dxN}
	    d\bsx^{[N]}_t =& ((\bA^{[N]} + \alpha_N \bbI) \bsx^N_t  \nonumber\\
        &+(\bB^{[N]}  + \beta_N \bbI) \bsu^{[N]}_t)dt + d\bsw^{[N]}_t,
	\end{align}
	where 
	\begin{equation}
        \bsw^N_t(\alpha) :=  \sum_{k=1}^N \bbS^N_k(\alpha) \sum_{r=1}^\infty c_{k,r}^N W_r(t),\label{C_defn}
    \end{equation}
    in which $\{W_r\}_{r=1}^\infty$ is a sequence of independent Brownian motions, and denote the $L_2[0,1]$ limit by 
    \begin{equation}
        \bsw^\infty_t(\alpha) := \lim_{N \to \infty}  \bsw^{[N]}_t(\alpha) = \sum_{r=1}^\infty \sqrt{\lambda_r} \phi_r(\alpha) W_r(t).
    \end{equation}
    
    By Mercer's Theorem (see, e.g. \cite{gohberg_basic_1981}), $\bQ$ has the eigenvalue and basis representation: 
    \begin{equation}
		\bQ(\alpha,\beta) t = \sum_{r=1}^\infty \sqrt{\lambda_r} \phi_r(\alpha)\phi_r(\beta) = \bE[\bsw^\infty_t(\alpha) \bsw^\infty_t(\beta)].
	\end{equation}

    To ensure the existence of this limit, we explicitly require the processes constructed in (\ref{C_defn}) to be Cauchy in the $L_2[0,1]$ norm, i.e., for all $\epsilon > 0$, there exists $M > N > N_0(\epsilon)$ such that
    \begin{align}
        &\bE[||\bsw^{[N]}_t - \bsw^{[M]}||^2_2]\\
        &\leq \sum_{r=1}^\infty \int_0^1\left(\sum_{j=1}^M \bbS^M_j(\alpha) c_{j,r}^M - \sum_{k=1}^N \bbS^N_k(\alpha) c_{k,r}^N\right)^2  d\alpha < \epsilon. \nonumber
    \end{align}

    Finally, we observe that as a result of the specifications above, the state process $\bsx^{[N]}_t$ on partition $P_i$  has a one-to-one correspondence with the state of the $i$th node of $x^N_t$ given by
    \begin{align}
        [x^N_t]_i:= \bsx^{[N]}_t(\alpha)  ~ \text{for} ~ \alpha \in P_i, ~ 1 \leq i \leq N.
    \end{align}

    \subsubsection{Common Graphons}
    There are a few common graphons that will be further investigated in Sec. V. In this subsection, a small set of common graphons and their associated dynamical systems properties are further investigated. In this paper, these are primarily used to generate random graphs using the W-random graph method \cite{lovasz_large_2012}.
    
    \textit{Erd\H{o}s-Renyi graphs} are one of the most common methods of generating random graphs. In an ER graph for every two vertices $i$ and $j$ in a graph of size $N$, an (undirected) edge $e_{i,j}$ exists with probability $p$. That is,
    \begin{align}
        \mathbb{P}(A^N_{i,j} = 1) = p, ~~ 1 \leq i,j \leq N.
    \end{align}
    From this, it is clear that the graphon limit of an ER graph is simply the constant function $\bW(\alpha,\beta) = p, ~ \alpha,\beta \in [0,1]$. 

    \textit{Uniform Attachment Graphs} \cite{lovasz_large_2012} are a more sophisticated increasing random graph model which possesses a smooth graphon limit. It is constructed inductively. Start with a single node graph $G^1$ (with associated adjacency matrix $A^1 = 0$). Then, given a UAG $G^{N-1}$, add a node and connect each pair of non-adjacent nodes with probability $\frac{1}{N}$ to create $G^N$.
    This has the graphon limit:
    \begin{align}
        \bW(\alpha,\beta) = 1 - \max(\alpha,\beta), ~~ \alpha,\beta \in [0,1].
    \end{align}
    \textit{Small World Graphs} \cite{watts_small_2018} model networks with a high level of local clustering, a low level of global clustering, and a low graph diameter. Medvedev \cite{medvedev_nonlinear_2014} presents one potential model of such graphs called a W-small world graphon. Here, we propose a limit model of small world graphs where the node connection probability is given by the sum of two truncated Gaussian functions of variance $\sigma^2$ on each horizontal in the unit square (\ref{eq:SWG_gaussian}), these are shifted by an offset $\gamma$ so that the resulting surface (\ref{eq:SWG_surface}) is symmetrically distributed with respect to the diagonal; it is normalized to have a maximum value of one on the diagonal,
    \begin{align}\label{eq:SWG_gaussian}
        \boldsymbol{G}_{SW}(\alpha,\beta) &=\exp \left(-\frac{1}{2}\bigg(\frac{\alpha - \beta}{\sigma}\bigg)^2 \right), ~~\alpha,\beta \in [0,1],\\
        \bW_{\text{SW}}(\alpha,\beta) &= 0.5 \boldsymbol{G}_{SW}(\alpha-\gamma,\beta) + 0.5\boldsymbol{G}_{SW}(\alpha,\beta-\gamma). \label{eq:SWG_surface}
    \end{align}
    Identifying 0 with the 0-degree position on a circle, 1 with the $\pi$ location and invoking symmetry shows this graphon shares some of the required SW network properties listed above.
    
    \begin{figure}
        \centering
        \includegraphics[width=\linewidth]{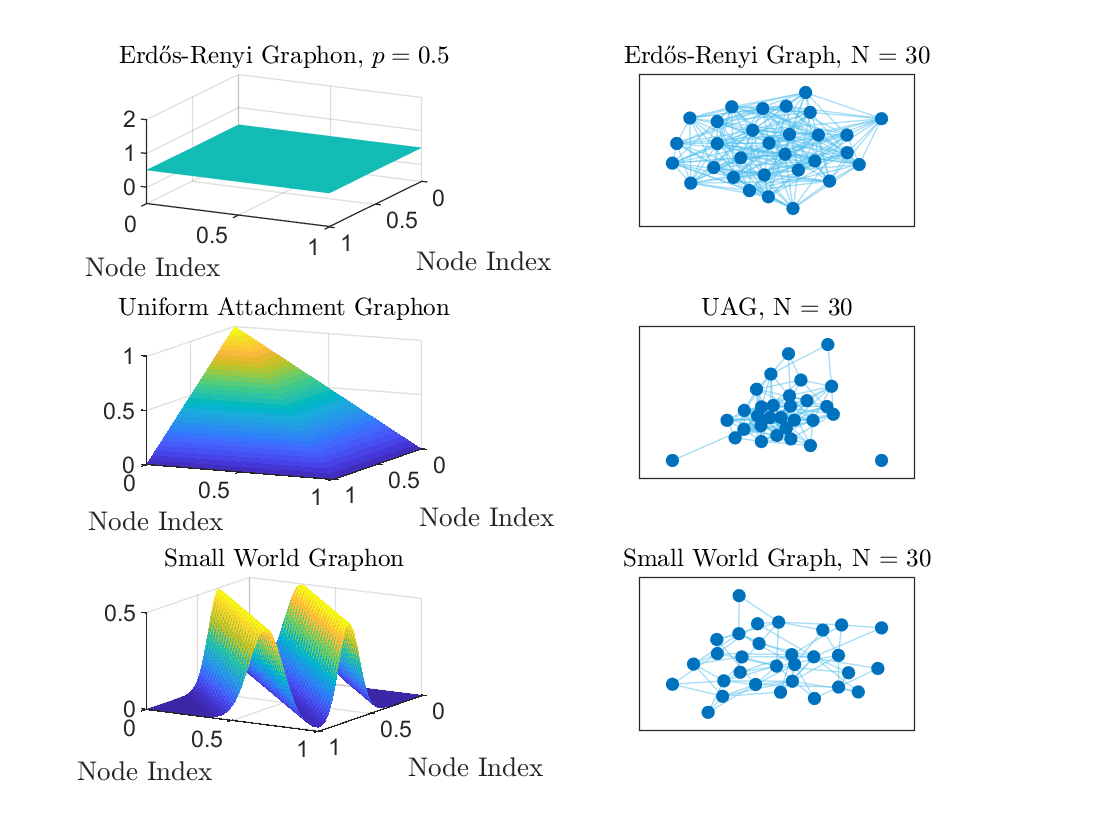}
        \caption{Top row: An Erd\H{o}s-Renyi graphon and sample graph. Middle row: A uniform attachment graphon and graph. Bottom row: A small world graphon.}
        \label{fig:enter-label}
    \end{figure}
    \textit{Low rank graphons} are a special class of graphons which possess a finite number of eigenfunctions. These are explored further in Sec. \ref{sec:low_rank}.
    
    The convergence of the finite order network system to the infinite limit graphon system can now be analyzed.
    \begin{theorem} \label{thm:mean_convergence}
        Let $\bsx^{[N]}_t$ solve the following graphon stochastic differential equation,
        \begin{equation}
            d\bsx^{[N]}_t = \bbA^{[N]} \bsx_t^{[N]}dt + d\bsw^{[N]}_t, ~~ \bsx^{[N]}_0 \in \mathcal{L}^2[0,1],
        \end{equation}
        where $\bQ^{[N]}$ is the covariance operator of $\bsw^{[N]}$. Let $\Phi^{[N]}(t,s)$ and $\Phi^{[M]}(t,s)$ refer to the semigroup solutions to the state systems $\bsx^{[N]}_t$ and $\bsx_t^{[M]}$ respectively. Assume for each triple $\epsilon_0,\epsilon_1, \epsilon_2 > 0$ there exists $N_0$ such that for all $N > M > N_0,$
        \begin{align}
            &||\bsx^{[N]}_0 - \bsx^{[M]}_0||_2^2 < \epsilon_0, &\text{(A0)}\\
            &||\Phi^{[N]}(t,s) - \Phi^{[M]}(t,s)||^2_\op < \epsilon_1, &\text{(A1)}\\
            &\bE[||\bsw^{[N]}_t - \bsw^{[M]}_t||^2_2] < \epsilon_2, &\text{(A2)}
        \end{align}
        and there exist $\alpha, C < \infty,$  such that, for all $N,$ 
        \begin{align}
            &\sum_{k=1}^N \sum_{r=1}^N |c^N_{k,r}|^2 \leq C, &\text{(B0)}\\
            & \int_0^s||\Phi^{[N]}(t,s)||^2_\op ds \leq \alpha, ~~ t \in [0,T].  &\text{(B1)}
        \end{align}
        Then, for each $\epsilon >0$, there exists $ N_0'$ such that for all $N,M > N_0'$ 
        \begin{equation}
            \bE[||\bsx^{[N]}_t - \bsx^{[M]}_t||_2] < \epsilon.
        \end{equation}
        Hence, there exists an ${L}^2[0,1]$  limit process $\bsx^\infty_t$ constituting the unique mild solution of (\ref{linear_sys}) satisfying $\bE[||\bsx_t^{[N]} - \bsx^\infty_t||_2] \to 0$ as $N$ goes to infinity, that is, 
        \begin{align}
            d\bsx_t = \bbA \bsx_t + d\bsw_t, ~ \bsx_0 \in L_2[0,1].
        \end{align}$\hfill \square$
    \end{theorem}
    The proof relies on the convergence of the operator norms of $\Phi^{[N]}(t,s)$ to $\Phi^{[N]}(t,s)$ implying $L_2$ convergence of $\Phi^{[N]}(t,s) \boldsymbol{v}$ for any $\boldsymbol{v} \in L_2[0,1]$ and is given in Appendix \ref{app:mean_convergence}.
    
    \textit{Remark:} This result was provided in \cite{dunyak_linear_2022} in the time-invariant operator case. When $\bbA^{[N]}_t = \bbA^{[N]}$ and $\bbA^{[M]}_t = \bbA^{[M]}$ for all $t \in [0,T]$, assumption (A1) can be relaxed to
    \begin{align}
        ||\bbA^{[N]} - \bbA^{[M]}||^2_\op < \epsilon_1,
    \end{align}
    i.e., the piecewise constant graphons corresponding to the finite graphs $A^N$ and $A^M$ converge in the $L_2$ operator norm.



\section{Linear Quadratic Control}
    \subsection{Finite Time Horizon}
    
    Given a linear stochastic dynamical system of type (\ref{linear_sys_B}), a linear quadratic Gaussian optimal control problem with Q-noise is given by the following performance function:
    \begin{align} \label{eq:objective_function}
        J(\bsu,\bsx_0) =& \bE[\int_0^T  (\bsx^*_t \bbM \bsx_t + \bsu^*_t \bbR \bsu_t) dt + \bsx^*_T \bbM_T \bsx_t], \nonumber\\
        &~~\bsx_t,\bsu_t \in L^2[0,1],
    \end{align}
    where $\bbM= \bM + m \bbI$ and $\bbR = \bR + r \bbI$ are bounded functionals composed of a non-negative compact operator $\bM \geq 0$ and $\bR > 0$, with $m \geq 0$ and $r > 0$ respectively.
    A control input $\bsu_t$ is admissible when it is adapted to the sigma algebra generated by $\bsw_t$ and where $\int_0^T |\bsu_t|^2dt < \infty$. Problems of this variety will henceforth be referred to as Q-LQG problems, and denoted by their system parameters $\{\bbA,\bbB,\bQ,\bbM,\bbM_T,\bbR\}$.

    \begin{theorem}
        Suppose that $\bbM,\bbR,\bbM_T$ are bounded positive self-adjoint $L^2[0,1]$ operators, and that $\bbR: L_2[0,1] \to L_2[0,1]$ is invertible.
        Then, the performance function (\ref{eq:objective_function}) is minimized with $\bsu_t = -\bbR^{-1} \bbB^* S_t \bsx_t$, where $S:[0,1] \times [0,1] \times [0,T] \to \mathcal{R}$ is an $\mathcal{L}^2[0,1]$ linear operator for all $t$ whose kernel satisfies the following Riccati equation,
        \begin{align} \label{riccati}
            -\frac{d}{dt}\langle {S}_t \bsv,\bsv \rangle =& 2 \langle \bbA \bsv, S_t \bsv \rangle - \langle S_t\bbB \bbR^{-1} \bbB^* S_t \bsv, \bsv \rangle, \nonumber \\
            &+ \langle \bbM \bsv, \bsv \rangle, ~~~~\forall \bsv \in L_2[0,1] \nonumber \\
            S_T =& \bbM_T.
        \end{align}
    \end{theorem}
    \textit{Proof}: See \cite[Sec. 4]{ichikawa_dynamic_1979} with $F= \bbI, ~ D = 0, ~ C= 0$ (hence $\Gamma(\cdot) = 0$ and $\Delta(S_t) = 0$). $\hfill \square$

    \begin{corollary}[\cite{ichikawa_dynamic_1979}]
        Given a Q-LQG problem with parameters $\{\bbA,\bbB,\bQ,\bbM,\bbM_T,\bbR\}$ where $S_t$ is the self-adjoint $L_2[0,1]$ operator solving (\ref{riccati}), $S_t$ is unique in the space of $L^2[0,1]$ non-negative operators.
        
        In addition, the value function is given by
        \begin{align}
            V(\bsx_t,t) = \bsx_t^* S_t \bsx_t + \int_t^T \text{trace} (S_r \bQ) dr.
        \end{align}
        $\hfill \square $
    \end{corollary}

    As expected, the intensity of the Q-noise does not change the optimal control $\bsu_t$, but does impact the value function of the Q-LQG problem. In order to show that the finite dimensional linear quadratic Gaussian problem on a network converges to an infinite-dimensional Q-LQG problem in the sense of converging value functions, state, and control functions, it must first be shown that the solution $S_t^{[N]}$ of the piece-wise constant Q-LQG problem is bounded in operator norm uniformly in $N \in \mathcal{Z}$.

    \begin{lemma} \label{lemma:uniform_bounded_lemma}
        Let $A^N$ and $B^N$ be bounded self-adjoint $L^2[0,1]$ operators, $M^N$ and $M^N_T$ be bounded positive $L^2[0,1]$ operators, and $R^N$ be a bounded strictly positive $L^2[0,1]$ operator converging to $A,~B,~M,~M_T,$ and $R$ respectively in the operator norm sense for $\{A^N,B^N,M^N,M_T^N,R^N\}$. Let $S^N_t$ be a positive, self-adjoint $L^2$ operator satisfying
        \begin{align} \label{eqn:uniform_bounded_lemma}
            -\dot{S}^N_t &= A^N S^N_t + S^N_t A^N - S^N_t B^N R^{-1} B^N S^N_t + M^N  \nonumber\\
            S^N_T &= M^N_T.
        \end{align}
        Then, there exists $0 <c_N<\infty $ such that
        \begin{equation}
            ||S^N_t||_\op \leq 2||M_T||_\op +(T-t)c_N. 
        \end{equation}
    \end{lemma}
    \textit{Proof:} See Appendix \ref{app:uniform_bounded_lemma}.

    From Theorem 3.1, the minimizing controls to the limit Q-LQG problem and the piece-wise constant Q-LQG problem are respectively $\bsu_t = -\bbR^{-1} \bbB^* S_t \bsx_t$ and $\bsu_t^{[N]} = -{(\bbR^{[N]})}^{-1} \bbB^{{[N]*}} S^{[N]}_t \bsx^{[N]}_t$.
    
    \begin{theorem}[Q-LQG Convergence] \label{thm:LQG_convergence}
        Let $\bsx_t$ be a system of the form
        \begin{equation} \label{linear_sys_B}
            d \bsx_t = (\bbA \bsx_t + \bbB\bsu_t)dt + d\bsw_t,
        \end{equation}
        and let $\bsx^{[N]}_t$ be a system of the form
        \begin{equation}
            d \bsx_t^{[N]}= (\bbA^{[N]} \bsx_t^{[N]} + \bbB^{[N]}\bsu_t^{[N]})dt + d\bsw_t^{[N]},
        \end{equation}
        where $\bbA^{[N]} \to \bbA$, $\bbB^{[N]} \to \bbB$, and $\bQ^{[N]} \to \bQ$ in the $L^2$ operator norm sense, and let assumptions (A0)-(A2), (B0), and (B1) of Theorem \ref{thm:mean_convergence} be satisfied. In addition, let $\bbR^{[N]}$, $\bbM^{[N]}$, and $\bbM^{[N]}_T$ be bounded, positive, self-adjoint operators converging to $\bbR$, $\bbM$, and $\bbM_T$ in the operator norm sense. 
        
        Let $S_t$ and $S^N_t$ be the positive, bounded self-adjoint operators solving the functional Riccati equation \ref{riccati} for $(\bbA,\bbB,\bbR,\bbM,\bbM_T)$ and $(\bbA^{[N]},\bbB^{[N]},\bbR^{[N]},\bbM^{[N]},\bbM^{[N]}_T)$ respectively. 
        Then, for every $\epsilon >0$, there exists an $N(\epsilon)$ such that for all $N > N(\epsilon)$,
        \begin{equation}
            \bE[||\bsx_t^{[N]} - \bsx_t||_2] < \epsilon.
        \end{equation}
    \end{theorem}
    \textit{Proof:} See Appendix \ref{app:LQG_convergence}.
    
    As the control input and state trajectory converge in the $L_2$ sense for all time, the optimally controlled finite dimensional network performance function value
        \begin{align}
            \inf_{u^N_t} \bE\left[\int_0^T (x^N_t)^* M^N x^N_t + (u_t^N)^* R^N u_t^N dt\right]
        \end{align}
        converges to the value of the infinite dimensional graphon system value
        \begin{align}
            &\inf_{\bsu} \bE[\int_0^T  (\bsx^*_t \bbM \bsx_t + \bsu^*_t \bbR \bsu_t) dt + \bsx^*_T \bbM_T \bsx_t],\\
            &~~\bsx_0,\bsu_t \in L^2[0,1].\nonumber
        \end{align}

    \subsection{Long-Range Average}\label{sec:infinite_time_horizon}
    In contrast to deterministic systems, the infinite time horizon optimal control problem does not have a finite value. Hence, we consider the long-range average Q-LQG problem given by:
    \begin{align}
        V_\infty(\bsx_0) :=&  \inf_{\bsu} \lim_{T \to \infty} J_T(\bsu,\bsx_0) \\=
        &\inf_{\bsu} \lim_{T \to \infty} \frac{1}{T}\bE[\int_0^T  (\bsx^*_t \bbM \bsx_t + \bsu^*_t \bbR \bsu_t) dt], \nonumber\\
            =& \lim_{T \to \infty} \frac{1}{T} ( \bsx_0^* S_t \bsx_0 + \int_0^T \text{trace} (S_r \bQ) dr). \nonumber
    \end{align}
    The solution to the long-range average Q-LQG problem is given via the unique positive solution $S_\infty$ to the algebraic Riccati equation:
    \begin{align}
        0 =& 2 \langle \bbA \bsv, S_\infty \bsv \rangle - \langle S_\infty\bbB \bbR^{-1} \bbB^* S_\infty \bsv, \bsv \rangle\\
            &+ \langle \bbM \bsv, \bsv \rangle, ~~~~\forall \bsv \in L_2[0,1]. \nonumber
    \end{align}
    The solution $S_t$ to equation (\ref{riccati}) converges to $S_\infty$ exponentially, yielding
    \begin{align}
        V_\infty(\bsx_0) = V_\infty := \text{trace} (S_\infty \bQ) ~~ \forall \bsx_0 \in L_2[0,1].
    \end{align}
    In particular, when $\bbB = \bbR = \bbI$, and $\bbA$ is symmetric, this can be solved with
    \begin{align} \label{eq:alg_riccati_soln}
        S_\infty = \bbA +(\bbA^2 + \bbM)^{\frac{1}{2}},
    \end{align}
    where, for a positive operator $M$, $(M)^{\frac{1}{2}}$ is the positive operator solving $(M)^{\frac{1}{2}*}(M)^{\frac{1}{2}} = M$.

    Note that the rate of loss is given by $\text{trace}(S_\infty \bQ)$ which is the Hilbert-Schmidt inner product. Given an orthornormal basis of $L_2[01]$, $\{\phi_k\}_k^\infty$:
    \begin{align}
        \text{trace}(S_\infty \bQ) = \sum_k \langle S_\infty \phi_k, \bQ \phi_k \rangle.
    \end{align}

    Restricting ourselves to the case where $\bbM = \bbI$, $||\bQ||_{\text{HS}} = 1$ and $||\bA||_{\text{HS}}  = 1$, this is maximized for the eigenfunction corresponding to the largest eigenvalue of $\bA$.

    \begin{lemma}[Maximum Trace Lemma]\label{lemma:lta_trace}
        Let $\bbM = \bbI$, and let $\{\phi_k\}_{k=0}^\infty$ be the set of orthonormal eigenvectors of $\bA$ with eigenvalues $\{\lambda_k\}_{k=0}^\infty$. Let $\overline{\lambda} = \sup_k \lambda_k$ and $\underline{\lambda} = \inf_k \lambda_k$ with associated eigenfunctions $\overline{\phi}$ and $\underline{\phi}$ respectively if $\overline{\lambda}$ and $\underline{\lambda}$ are obtained for a finite $k$. Then, for systems driven by of the form $\bbA = \bA + a \bbI$,
        \begin{align}
            \sup_{||\bQ||_{\text{HS}} = 1} \text{\emph {trace}}(S_\infty \bQ) =& \sup_{||\bQ||_{\text{HS}} = 1}\sum_k \langle S_\infty \phi_k, \bQ \phi_k \rangle \\
            =& (\overline{\lambda} + a) + \sqrt{(\overline{\lambda} + a)^2 +1}. \label{eq:trace_value}
        \end{align}
        and attains the supremum for $\bQ = \langle \cdot,\overline{\phi} \rangle$ when $\overline{\lambda}$ is obtained for a finite $k$, and respectively obtains the infimum for $\bQ = \langle \cdot,\underline{\phi} \rangle$ when $\underline{\phi}$ is attained for finite $k$.
    \end{lemma}
    \textit{Proof:} Note that the operator $(\bbA^2 + \bbI)$ can be expressed as
    \begin{align}
        &\bbA^2 + \bbI = \bA^2 + 2 a \bA + (a^2 + 1) \bbI\\
        &= \sum_{k = 1}^\infty \lambda_k^2 \langle \cdot, \phi_k \rangle \phi_k + 2 a \lambda_k \langle \cdot, \phi_k \rangle \phi_k + (a^2 + 1)  \langle \cdot, \phi_k \rangle \phi_k \nonumber \\
        &= \sum_{k=1}^\infty ( \lambda_k^2 + 2a\lambda_k + a^2 + 1) \langle \cdot, \phi_k \rangle \phi_k \nonumber\\
        & = \sum_{k=1}^\infty ( (\lambda_k + a)^2 + 1) \langle \cdot, \phi_k \rangle \phi_k. \nonumber
    \end{align}
    Taking the positive root of the (necessarily positive) eigenvalues gives us the operator root,
    \begin{align}
        (\bbA^2 + \bbI)^{\frac{1}{2}} = \sum_{k=0}^\infty \sqrt{(\lambda_k + a)^2+ 1} \langle \cdot, \phi_k\rangle,
    \end{align}
    which then yields
    \begin{align}
        S_\infty = \sum_{k=0}^\infty \left((\lambda_k + a) + \sqrt{(\lambda_k + a)^2 + 1} \right) \langle \cdot, \phi_k\rangle \phi_k.
    \end{align}
    First, suppose that $\overline{\phi}$ is attained for a finite $k$. Then, setting $\bQ = \langle \cdot,\overline{\phi} \rangle$ yields
    \begin{align}
        \text{trace}(S_\infty \bQ) = (\overline{\lambda} + a)+ \sqrt{(\overline{\lambda} + a)^2 + 1}.
    \end{align}
    As $\bQ$ is constrained by $||\bQ||_{\text{HS}} = 1$ and the eigenbasis is orthonormal, assigning any positive value to a different eigenfunction cannot increase this value.

    Likewise, if the infimum $\underline{\lambda}$ is attained for finite $k$, the minimum value is attained for $\bQ = \langle \cdot, \underline{\phi}\rangle $ and
    \begin{align}
        \text{trace}(S_\infty \bQ) = (\underline{\lambda} + a)+ \sqrt{(\underline{\lambda} + a)^2 + 1}.
    \end{align}
    If the supremum (or infimum, respectively) is not attained for finite $k$, then the limit $\lambda_{k \to \infty} = 0$ implies that the supremum (respectively infimum) value is
    \begin{align}
        \text{trace}(S_\infty \bQ) = a+ \sqrt{ a^2 + 1}
    \end{align}
    which can be made arbitrarily close by setting $\bQ = \langle \cdot, \phi_k \rangle$ for arbitrarily large $k$.
    $\hfill \square$
    
    \subsection{Exponential Discounting}
    The infinite time horizon problem with exponential discounting is presented below. Unlike the long-term averaging problem, the value function is found by explicitly solving a Hamilton-Jacobi-Bellman equation.
    \begin{lemma}[Infinite Horizon Discounting]\label{lemma:discount_trace}
        Let $\rho > 0$. Then the infinite horizon discounted cost performance functional is given by
        \begin{align}
            J_\rho(\bsx_0, \bsu) = \bE[\int_0^\infty e^{-\rho t} (\bsx_t^* \bbM \bsx_t + \bsu_t^* \bbR \bsu_t)dt].
        \end{align}
        This is minimized by $\bsu_t = -\bbR^{-1}\bbB^* S_\infty$, where $S_\infty$ solves the discounted Algebraic Riccati equation
        \begin{align}
            \rho S_\infty = S_\infty \bbA^* + \bbA S_\infty - S_\infty \bbB \bbR^{-1} \bbB^* S_\infty + \bbM.
        \end{align}
        When $\bbB = \bbR = \bbI$ and $\bbA$ is symmetric, the unique positive symmetric solution is given by
        \begin{align}
            S_\infty^\rho = (\bbA - \frac{ \rho }{2} \bbI) +((\bbA - \frac{ \rho }{2}\bbI)^2 + \bbM)^{\frac{1}{2}}.
        \end{align}
    \end{lemma}
    \textit{Proof:} For the existence and uniqueness of the solution, refer to \cite[Section 2.6.1.3]{fabbri_stochastic_2017}. For a function $V:L_2[0,1] \to \calR^1$, let $DV$ be the Frechet derivative of $V$. Define the value function $V_\rho:L_2[0,1] \to \mathcal{R}^1$
    \begin{align}
        V_\rho(\bsx) = \langle \bsx, S^\rho_\infty \bsx \rangle + \frac{1}{ \rho} \text{trace}(\bQ S^\rho_\infty), ~~\bsx \in L_2[0,1]
    \end{align}
    be a classical solution of the Hamilton-Jacobi-Bellman equation of an infinite horizon discounted cost performance problem, satisfying
    \begin{align}
        \rho V_\rho &- \frac{1}{2}\text{trace}(\bQ D^2 V_\rho) - \langle \bbA \bsx, DV_\rho \rangle \\
        &- \inf_{\bsu} \{  \langle \bbB \bsu, DV_\rho \rangle + \langle \bsx, \bbM \bsx \rangle + \langle \bsu, \bbR \bsu\rangle  \} = 0. \nonumber
    \end{align}
    After evaluating the Frechet derivatives $DV_\rho = 2 S_\infty^\rho \bsx$ and $D^2 V_\rho = 2S_\infty^\rho$, this is equivalent to
    \begin{align}
        \rho &(\langle \bsx, S^\rho_\infty \bsx \rangle + \frac{1}{\rho } \text{trace}(\bQ S^\rho_\infty)) \\
        &- \text{trace}(\bQ S^\rho_\infty) - 2 \langle \bbA \bsx, S^\rho_\infty \bsx \rangle \nonumber \\
        &- \inf_{\bsu} \{ 2 \langle \bbB \bsu, S^\rho_\infty \bsx \rangle + \langle \bsx, \bbM \bsx \rangle + \langle \bsu, \bbR \bsu\rangle  \} = 0. \nonumber
    \end{align}
    Noting that the infimization holds for $\bsu = - \bbR \bbB^* S^\rho_\infty \bsx$, this is equivalent to
    \begin{align}
         \langle \bsx,  \rho S^\rho_\infty \bsx \rangle =& 2 \langle \bbA \bsx, S^\rho_\infty \bsx \rangle - \langle \bbB \bbR^{-1} \bbB S^\rho_\infty \bsx, S^\rho_\infty \bsx \rangle,  \nonumber \\
         & + \langle \bsx, \bbM \bsx \rangle, ~~ \forall \bsx \in L_2[0,1]
    \end{align}
    and the result follows. $\hfill \square$

    \textit{Remark:} In the special case of $\bbB = \bbR = \bbM = \bbI$, the optimal discounted performance control decreases the control input along all the direction of all eigenfunctions. This results in a weaker control input for all actuators compared to the long-range averaging solution.

    To illustrate Lemma \ref{lemma:lta_trace} and Lemma \ref{lemma:discount_trace}, a selection of worst case scenarios are presented in Table 1, comparing the $S_\infty$ value to the $S_t$ value at $T = 1$ (with the local forcing term $a = 0$). In each example, the underlying graphon $\bA$ has only non-negative eigenvalues, and the best case scenario has the same cost (namely, $\text{trace}(S_\infty \bQ) = 1$). Notably, the Erodos-Renyi Graphon  $\bA(\alpha,\beta) = 0.5$ and $\bA(\alpha,\beta) = \cos(2\pi (\alpha-\beta))$ have the same worst case value, due to the fact that they have the same maximum eigenvalue.
    \begin{table*}
        \begin{center}
            \begin{tabular}{|c|c|c|c|}
                \hline
                Graphon &  Max eigenvalue of $\bA$ & $||S_\infty \bQ||^2_{\text{H.S.}}$ &$ ||S^\rho_\infty \bQ||_{\text{H.S.}}$, $\rho = 1$\\
                \hline \hline
                $\bA(x,y) = (x^2  - 1)(y^2 -1)$ & 0.533 & 1.666 & 1.034\\
                \hline
                $\bA(x,y)  = 0.5$ (Erdos-Renyi) &  0.5 & 1.618 & 1.000\\
                \hline
                $\bA(x,y)  = \cos(2\pi (x-y))$ & 0.5 & 1.618 & 1.000\\
                \hline
                $\bA(x,y)  = 1 - \max(x,y)$ (UAG) & 0.405 & 1.484 & 0.910\\
                \hline
                S.W. ($\sigma = 0.1, \gamma = 0.3$) & 0.183 & 1.200 & 0.783\\
                \hline
            \end{tabular}
        \end{center}
        \caption{A comparison of the worst case performance of various graphons with Hilbert-Schmidt norm bounded noise covariance $\bQ$. Calculating the H.S. norm $(\text{trace}(S_\infty \bQ))^2$ agrees with the maximum value calculated by Eq. (\ref{eq:trace_value}). As expected, the discounting problem with discount factor $\gamma =1$ has a lower expected cost than the long-time averaging problem.}
    \end{table*}

    When the relevant operators are infinite dimensional, the system cannot be fully simulated. However, when the operators are finite dimensional, the system can be fully analyzed as an $N$-dimensional system.
   
\section{Low Rank Graphons}
\label{sec:low_rank}
Here, the theory described in \cite{gao_subspace_2021} is extended to Q-noise systems. Define an invariant subspace of a linear operator $\mathbb{T}$ by $\calS \subset L_2[0,1]$ to be a subspace of $L_2[0,1]$ such that, for all $\bsx \in \calS$,
\begin{equation}
    \bsx \in \calS \implies \mathbb{T} \bsx \in \calS.
\end{equation}
Denote the orthogonal complement of $\mathcal{S}$ to be the subspace $\calS^\perp$, such that for all $\boldsymbol{v} \in \calS$ and all $\Breve{\boldsymbol{v}} \in \calS^\perp$, $\langle  \boldsymbol{v}, \Breve{\boldsymbol{v}}\rangle = 0$. This partitions $L_2[0,1]$ into the two orthogonal spaces, $\calS$ and $\calS^\perp$, and hence $L_2[0,1] = \calS \oplus \calS^\perp$. An operator $\mathbb{T}$ is said to be low rank (with respect to an invariant subspace $\calS$) when, for all $\Breve{\boldsymbol{v}} \in \calS^\perp$, $\mathbb{T} \Breve{\boldsymbol{v}} = 0$.

For a given linear Q-noise graphon process as generated by (\ref{linear_sys}), with $\bbA = \bA + a \bbI$, $\bbB = \bB + b\bbI$, and a Q-noise $\bsw_t$ with covariance operator $\bQ$, make the following assumptions:

\begin{itemize}
    \item (\textbf{LRG1}) $\calS$ is spanned by a finite number of orthonormal $L_2[0,1]$ functions denoted $\bsf = (\bsf_1,... \bsf_N)$. 
    \item (\textbf{LRG2}) The operators $\bA$ and $\bB$ are finite rank self-adjoint graphon operators which share the non-trivial invariant subspace $\calS$. That is, for all $\Breve{\boldsymbol{v}} \in \calS^\perp$, $\bA \Breve{\boldsymbol{v}} = 0$ and $\bB \Breve{\boldsymbol{v}} = 0$.
    \item \textbf{(LRG3)} $\bsw_t$ is finite dimensional, and has a representation of the form
    \begin{align}
        \bsw_t = \sum_{k=1}^N \sqrt{\lambda_k} \bsf_k W^k_t.
    \end{align}
    Equivalently, the covariance operator $\bQ$ is low rank with respect to $\mathcal{S}$.
\end{itemize}

By (LRG1), the state process of a linear Q-noise graphon process as generated by (\ref{linear_sys}) can be decomposed into two orthogonal components
\begin{align}
    \bsx_t =: \bsxf_t + \bsxp_t, ~~ \bsxf \in \calS, ~~ \bsxp \in \calS^\perp
\end{align}
where $\bsxf_t = (\bsx_t | \calS)$ is the orthogonal projection of $\bsx_t$ into $\calS$ and $\bsxp_t = (\bsx_t | \calS^\perp)$ is the orthogonal projection into $\calS^\perp$. Hence, $\bsxf_t$ consists of a linear combination of elements of $\calS$ and $\bsxp$ consists of a linear combination of elements of $\calS^\perp$. Similarly, decompose $\bsu_t$ into
\begin{align}
    \bsu_t =: \bsuf_t + \bsup_t.
\end{align}
By the Q-noise axioms, $\bsw_t$ is defined as a sum of weighted Wiener processes, each associated with an eigenbasis of $\bQ$. Consequently, $\bsw_t$ can also be decomposed into its orthogonal components with respect to $\calS$ and $\calS^\perp$. This has the general form
\begin{align}
    \bsw_t =& \bsw^\bsf_t + \Breve{\bsw}_t,
\end{align}
where
\begin{align}
    \bswf_t :=& (\bsw_t | \calS)\\
    =&\sum_{r=1}^N \langle \sum_{k=1}^\infty \sqrt{\lambda_k}\phi_k W^k_t, \bsf_r \rangle \bsf_r \nonumber\\
        =&\sum_{r=1}^N \bigg(\sum_{k=1}^\infty \sqrt{\lambda_k}\langle \phi_k, \bsf_r \rangle W^k_t\bigg) \bsf_r, \nonumber \\
    \bswp_t :=& (\bsw_t | \calS^\perp) = \bsw_t - \bswf_t. 
\end{align}
By (LRG3), $\bsw_t$ is low rank with respect to the common invariant subspace $\calS$ of $\bA$ and $\bB$, hence $\bswf_t = \bsw$, $\bswp_t = 0$, and $\bsxp_t$ is deterministic. Consequently, these processes evolve according to
\begin{align}
    d\bsxf_t &= ((\bA + a\bbI ) \bsxf_t + (\bB + b\bbI )\bsuf_t )dt + d\bsw_t,\\
    d\bsxp_t &= (a \bsxp_t  + b \bsup_t)dt,\\
    \bsxf_0 &\in \calR^N, ~~ \bsxp_0 \in L_2[0,1].
\end{align}
Notably, by the low rank assumptions on $\bA$ and $\bB$, the orthogonal process $\bsxp_t$ is diagonal---each point on the unit interval evolves as a single-dimensional linear differential equation. Because of this decomposition, the system can be modelled as a finite dimensional system.
\subsection{Projections onto the Invariant Subspace $\calS$}
To project the low rank linear Q-noise graphon system to the finite dimensional invariant subspace $\mathcal{S}$, define the $\mathcal{R}^N$-valued state processes $x^{\bsf}_t$, $u^\bsf_t$:
\begin{align}
    x^\bsf_t &:= [\langle \bsxf_t, \bsf_1\rangle, \langle \bsxf_t, \bsf_2 \rangle, ..., \langle \bsf_N, \bsxf_t\rangle],\\
    u^\bsf_t &:= [\langle \bsuf_t, \bsf_1\rangle, \langle \bsuf_t, \bsf_2\rangle, ..., \langle \bsf_N, \bsuf_t \rangle],
\end{align}
i.e., $x^\bsf_t$ and $u^\bsf_t$ are projections onto the coordinate space defined by $\bsf$.

Similarly, define the following $N \times N$ matrices 
\begin{align}
    &A_{ij} := \langle \bA\bsf_i, \bsf_j \rangle, ~~ B_{ij} := \langle \bB\bsf_i,  \bsf_j \rangle,\\
    &Q = \text{diag}(\{\lambda_k\}_{k=1}^N),
\end{align}
and let $W^\bsf_t$ be an $N$ dimensional Wiener process with covariance matrix $Q$. Then the state process $x^\bsf_t$ equivalently evolves according to the finite dimensional differential equation
\begin{align}
    dx^\bsf_t &= ((A + aI) x^\bsf_t + (B + bI)u^\bsf_t)dt + dW_t^\bsf,\\
    x^\bsf_0 &=[\langle \bsxf_0, \bsf_1\rangle, \langle \bsxf_0, \bsf_2 \rangle, ..., \langle \bsf_N, \bsxf_0 \rangle].
\end{align}

This construction allows for a low-dimensional analysis of $\bsxf_t$ and $\bsuf_t$, which can then be mapped back into the $L_2[0,1]$ space by associating each element $[u_t^\bsf]_k$ with its respective basis function $\bsf_k$,
\begin{align}
    \bsuf_t = \sum_{k=1} [u^\bsf_t]_k \bsf_k =: u^\bsf_t \circ \bsf.
\end{align}

This approach is particularly useful for linear quadratic problems on low rank graphon systems.

\subsection{Low Rank Linear Quadratic Control}
Consider a graphon Q-LQG optimal control problem with the operators $\bbA = \bA + a \bbI$ and $\bbB = \bB + b\bbI$, and cost operators $\bbM$, $\bbM_T$, and $\bR$.

As with the standard Q-LQG problem, the objective function is
\begin{align}
    J(\bsu,\bsx_0) = \bE[\int_0^T (\bsx_t^* \bbM \bsx_t + \bsu_t^* \bR \bsu_t) dt + \bsx_T^* \bbM_T \bsx_T].
\end{align}
By taking the orthogonal decomposition of $\bsx_t$ and $\bsu_t$,
\begin{align}
    \bsx_t^* \bbM \bsx_t =& (\bsxf_t + \bsxp_t)^* \bbM (\bsxf_t + \bsxp_t)\\
    =& (\bsxf_t)^* \bbM \bsxf_t + (\bsxp_t)^* \bbM \bsxp_t \nonumber 
\end{align}
and
\begin{align}
    \bsu_t^* \bR \bsu_t  = (\bsuf_t)^* \bR (\bsuf_t) +  (\bsup_t)^* \bR (\bsup_t).
\end{align}

As with the $\bA$ and $\bB$ operators, the $N \times N$ real matrices $M$ and $R$ can be defined as
\begin{align}
    &M_{ij} := \langle \bbM\bsf_i, \bsf_j \rangle, ~~ R_{ij} := \langle \bR\bsf_i,  \bsf_j \rangle,\\
    &{M_T}_{ij} = \langle \bbM_T\bsf_i, \bsf_j \rangle.
\end{align}
Then, the optimal control problem can be decomposed into the following $N$-dimensional LQG optimal control problem (which can be solved using standard Riccati equation methods) and an $L_2[0,1]$ deterministic orthogonal process,
\begin{align}
    J(\bsu, \bsx_0) =& J^\bsf(u^\bsf,x^\bsf_0) + \Breve{J}(\bsup, \bsxp_0),\\
    J^\bsf(u^\bsf,x^\bsf_0) :=& \bE[\int_0^T ((x^\bsf_t)^* M x^\bsf_t + (u^\bsf_t)^* R u^\bsf_t) dt \\
    &+ (x^\bsf_T)^* M_T x^\bsf_T], \nonumber \\
    dx^\bsf_t =& ((A + aI) x^\bsf_t + (B + bI)u^\bsf_t)dt + dC_t^\bsf,\\
    x^\bsf_0 =&[\langle \bsxf_0, \bsf_1\rangle, ..., \langle \bsxf_0, \bsf_N \rangle], ~~ x^\bsf_0 \in \calR^N,\\
    \Breve{J}(\bsup, \bsxp_0) =& \int_0^T (\bsxp_t^* \bbM \bsxp_t + \bsup_t^* \bR \bsup_t) dt + \bsxp_T^* \bbM_T \bsxp_T,\\
    d\bsxp_t =& (a \bsxp_t  + b \bsup_t)dt.\\
    \bsxp_0 \in& L_2[0,1]
\end{align}
Further, when $\bbM$, $\bbM_T$, and $\bR$ are low rank with respect to $\bsf$ except for a diagonal constant, then the minimizing solution to $\Breve{J}(\bsup, \bsxp_0)$ is effectively one-dimensional, as the feedback control is diagonal with identical coefficients for each $\alpha \in [0,1]$.
\begin{theorem}
    \label{thm:finite_rank_lqg}
    Define a Q-LQG problem with coefficients $\{\bA + a\bbI, \bB + b\bbI, \bQ, \bM + m\bbI, \bM_T + m_t \bbI, \bR + r\bbI\}$ where all operators are rank $N$ with respect to an orthogonal subspace $\mathcal{S}$, whose projections onto $\calS$ are denoted $\{A,B,Q,M,M_T,R\}$ respectively. Then, the optimizing control $\bsu_t^0$ is given by 
    \begin{align}
        \bsu_t^0 :=& \bsuf_t + \bsup_t,~~ \bsu_t^0 \in L_2[0,1], ~~ t \in [0,T]\\\
        \bsuf_t :=& \sum_{k=1} [-(R + rI)^{-1} (B + bI)^* P_t x^\bsf_t]_k \bsf_k ,\\
        \bsup_t:=& -\frac{b^2}{r} p_t \bsxp_t,
    \end{align}
    where $P_t$ and $p_t$ are time-varying operators solving the following $N$-dimensional and one dimensional Riccati equations respectively,
    \begin{align} \label{thm:finite_rank_riccati}
        -\dot{P_t} =& (A+aI)^*P_t + P_t(A+aI)\\
        &- P_t(B+bI)^*(R+rI)^{-1}(B+bI)P_t \nonumber \\
        &+ (M+mI), \nonumber \\
        P_T =& (M_T + m_T I), ~~ P_t \in \mathcal{R}^{N \times N}\label{thm:finite_rank_riccati_end}\\
        \dot{p}_t =& 2ap_t - \frac{b^2}{r}  p_t^2 + m,\label{thm:perpendicular_riccati}\\
        p_T =& m_T. \label{thm:perpendicular_riccati_end}
    \end{align}
\end{theorem}
\textit{Proof:} This result is analogous to the deterministic finite dimensional graphon LQR solution shown in \cite{gao_subspace_2021}. The Riccati equations (\ref{thm:finite_rank_riccati}-\ref{thm:finite_rank_riccati_end}) give the standard $N$-dimensional and one-dimensional operator solution $P_t$ and $p_t$ respectively. Using the solution $P_t$, the optimizing controls for the finite-rank subspace LQG problem (with respect to the $N$-dimensional state vector $x^\bsf_t$ is given by
\begin{align}
    u^\bsf_t =& -(R + rI)^{-1} (B + bI)^* P_t x^\bsf_t.
\end{align}
By associating the $k$th entry of the control vector with the corresponding basis function $\bsf_k$, the finite dimensional controls can be mapped back to the original space.

Similarly, while the orthogonal state complement $\bsxp_t$ is infinite dimensional, the feedback gain is one-dimensional due to the diagonal nature of the state process , and the optimal control can be found with by solving the scalar Riccati equations (\ref{thm:perpendicular_riccati}-\ref{thm:perpendicular_riccati_end}). 

This gives the optimal controls for $\bsuf_t$ and $\bsup_t$, and hence for $\bsu_t^0$. $\hfill\square$

\section{Numerical Examples}

For the following numerical simulations, the unit interval $[0,1]$ is partitioned into $N$ segments, and the $k$th partition segment $I_k$ is denoted as
\begin{align}
    I_1 := [0, \frac{1}{N}], ~~~ I_k := (\frac{k-1}{N},\frac{k}{N}].
\end{align}
In each example, $N=50$, and the state of the simulated systems follow the form
\begin{align}
    d\bsx^{[N]}_t = (\bbA^{[N]} \bsx^{[N]}_t + \bbB^{[N]} \bsu^{[N]}_t)dt + d\bsw^{[N]}_t,
\end{align}
as in equation (\ref{eq:dxN}). This discretized system is used as a approximate solution to the infinite dimensional system.

In the following sections, $\bbA = \bA + 0.1 \bbI$ and $\bbB = 0.1\bbI$, where $\bA$ is a symmetric graphon and $\bbI$ is the identity operator. To simulate the Q-LQG problems, we set a terminal time of $T=1$ and implemented Euler's method with a time increment of $\Delta t =0.001$.

There are three key results to be presented: first, the convergence of a linear quadratic Gaussian finite graph system to a graphon system. Next, that a graph system with a low rank graphon limit can be efficiently represented by a low rank decomposition. Finally, we demonstrate that the finite time horizon feedback solution converges to the infinite time horizon solution. In order to compare trajectories, we introduce the root squared distance of two system trajectories $\bsx_t$ and $\bsy_t$ at time $t$,
\begin{align}
    rmd(\bsx_t, \bsy_t) = \sqrt{ \langle \bsx_t - \bsy_t,\bsx_t - \bsy_t \rangle }.
\end{align}

\subsection{Low rank finite graph convergence}
\begin{figure}
    \centering
    \includegraphics[width = \linewidth]{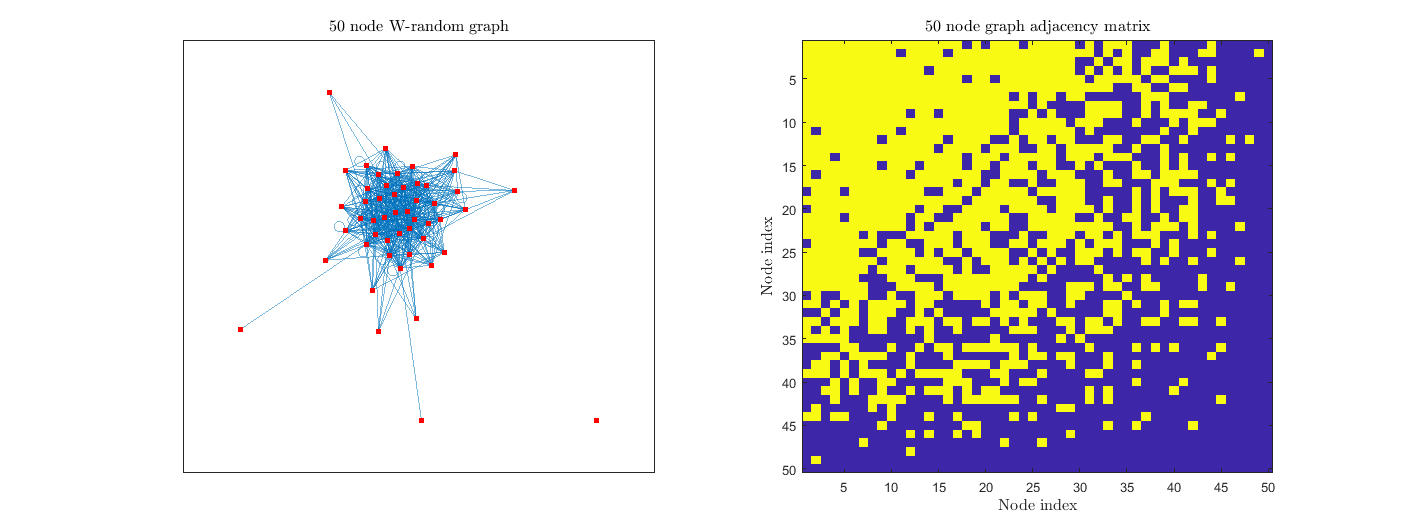}
    \caption{Left: a fifty node W-random graph. Right: the associated adjacency matrix to be used for the numerical simulations. Yellow squares represent an edge, blue squares represent a lack of an edge. The adjacency matrix is rank 49, despite the limit system being rank one.}
    \label{fig:50node_lowrank}
\end{figure}

Consider a finite graph generated using the following W-random graph \cite{lovasz_large_2012} kernel:
\begin{align}
    \bA(\alpha,\beta) = (\alpha^2 - 1)(\beta^2 - 1), ~~\alpha,\beta \in [0,1].
\end{align}
Clearly, $\bA$ is a rank one graphon, with a basis function given by 
\begin{align}
    \bsf(\alpha) = \frac{(\alpha^2 - 1)}{\sqrt{\int_0^1 (\beta^2 - 1)^2 d\beta}}, ~~ \alpha \in [0,1]
\end{align}
For each pair of gridpoints of the partition $I$, independently sample a Bernoulli random variable to generate an edge between that pair, with edge probability given by
\begin{align}
    \mathcal{P}(e_{ij} = 1) = \bA(\alpha_i, \alpha_j), ~ \alpha_i,\alpha_j \in [0,1], ~ i,j \in \{1,...,N\}.
\end{align} 
This creates a graph of 50 nodes, shown in Fig. \ref{fig:50node_lowrank} along with its adjacency matrix. Despite having a rank one limit, the finite adjacency matrix $A^{[N]}$ is full rank.

\begin{figure}
    \centering
    \includegraphics[width = 0.7\linewidth]{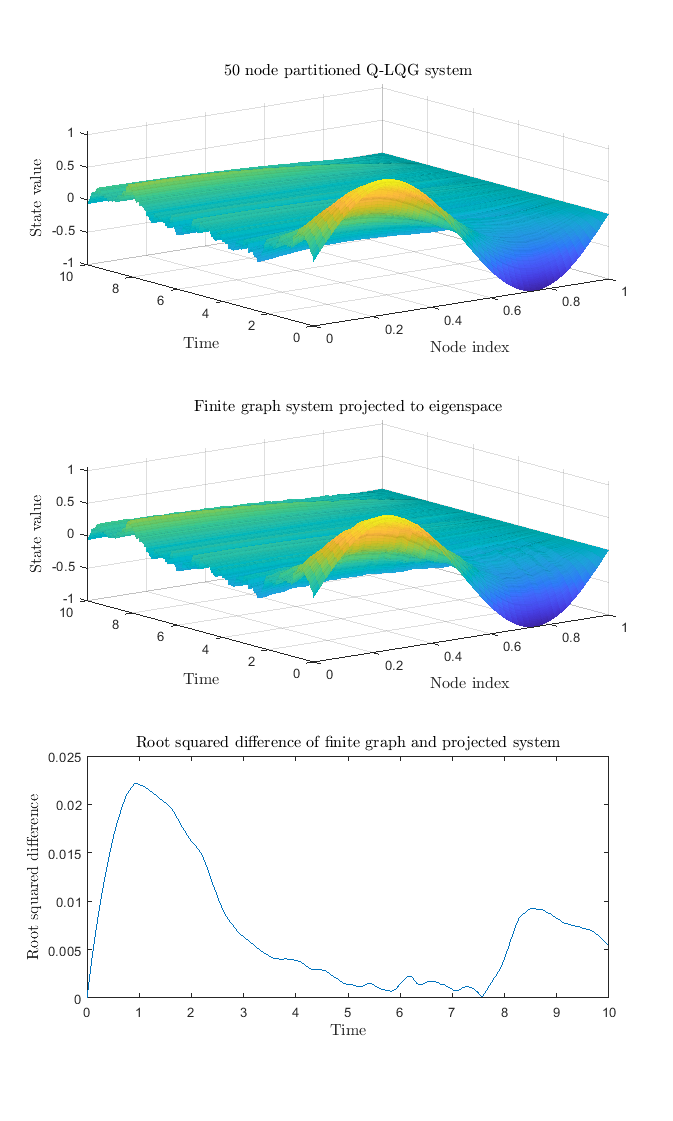}
    \caption{Top: the system generated with the finite graph using the piecewise constant graphon $\bA^{[N]}$. Middle: The system trajectory generated when $\bA^{[N]}$ is projected onto the eigenspace spanned by $\bsf$. Bottom: the root squared distance of the finite graph system trajectory and the projected graph system trajectory. The root squared distance has a maximum deviation of $0.023$, showing that the two trajectory surfaces are very similar.} 
    \label{fig:finite_graph_LQG}
\end{figure}

As with Section \ref{sec:infinite_time_horizon}, to simulate the worst case scenario, the Q-noise disturbance is placed on the basis function $\bsf$,
\begin{align}
    \bsw_t(\alpha) = \bsf(\alpha) W_t,  ~~ \alpha \in [0,1].
\end{align}
Then, the finite graph LQG problem can be solved with standard methods, and the limit system can be solved with Theorem \ref{thm:LQG_convergence}. The finite graph system trajectory is shown in Fig. \ref{fig:finite_graph_LQG}-I, and the low rank system created by projecting $\bA^{[N]}$ onto the normal basis function $\bsf$ is given by Fig. \ref{fig:finite_graph_LQG}-II. This is accomplished by simulating the system 
\begin{align}
    dx^\bsf_t &= \left((\langle \bA^{[N]} \bsf, \bsf \rangle + 0.1) x^\bsf_t + 0.1 u^\bsf_t\right)dt + dW^\bsf_t, \\
    d\bsxp_t &= (0.1 \bsxp_t + 0.1 \bsup_t)dt, ~~ t \in [0,1],\\
    & x^\bsf_0 \in \calR^1, ~ \bsxp_0 \in L_2[0,1],
\end{align}
where $\langle \bA^{[N]} \bsf, \bsf \rangle + 0.1$ is simply equal to the constant $1.7251$.
Even though the finite graph's adjacency matrix is full rank, the (piecewise constant) system projected onto the eigenspace spanned by $\bsf$ captures the behavior of the finite trajectory

\begin{figure}
    \centering
    \includegraphics[width = 0.7\linewidth]{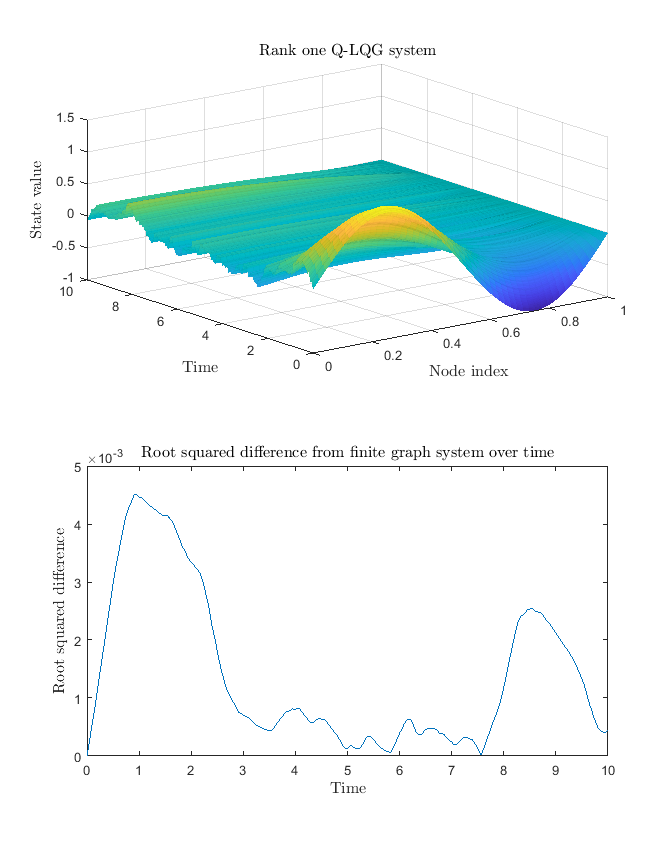}
    \caption{Top: the trajectory of the system under the rank one limit control. Bottom: the positive root of the squared distance between the finite graph system and the limit system over time.}
    \label{fig:limit_LQG}
\end{figure}

The limit system generated with $\bA$ in place of $\bA^{[N]}$ is shown in Fig. \ref{fig:limit_LQG}-I, and the root squared difference over time of the trajectory of the finite graph system and the limit system are shown in Fig. \ref{fig:limit_LQG}-II, calculated with

\subsection{Long-Range Average Comparison}
\begin{figure}
    \centering
    \includegraphics[width = 0.7\linewidth]{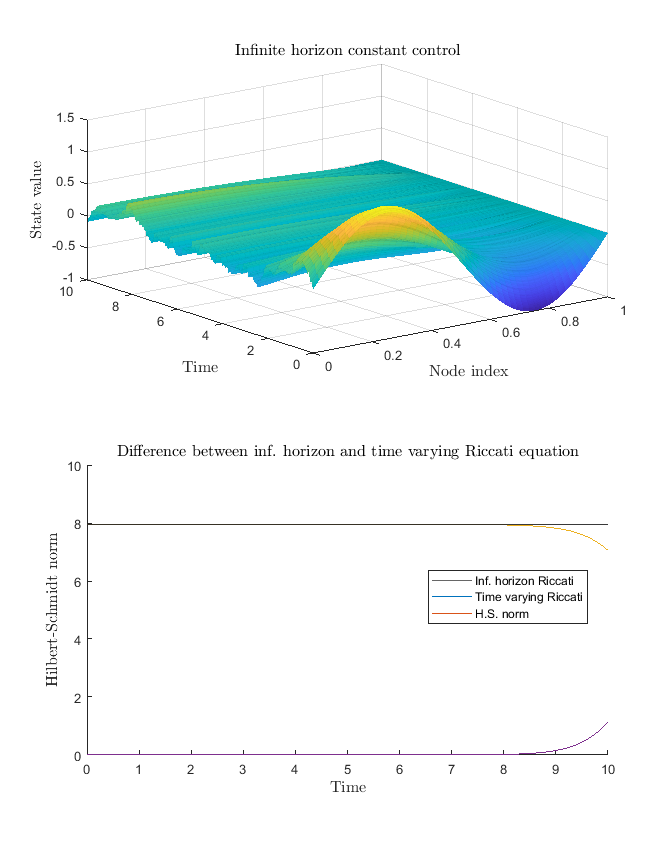}
    \caption{I: the trajectory of a system under infinite horizon control. II: The Hilbert-Schmidt norms of the infinite horizon Riccati equation solution and the time-varying Riccati equation solution. }
    \label{fig:inf_horizon}
\end{figure}

This can be found using the analysis of Section \ref{sec:infinite_time_horizon}. We apply the infinite horizon control found using the algebraic Riccati equation solution (\ref{eq:alg_riccati_soln}), and the resulting trajectory is shown for a terminal time of $T = 10$ in Fig. \ref{fig:inf_horizon}-I. A comparison of the Hilbert-Schmidt norms of both $S_t$, the time-varying solution to the differential Riccati equation associated with the system and the infinite horizon solution of the algebraic Riccati equation is shown in Fig. \ref{fig:inf_horizon}-II. The time varying Riccati solution converges exponentially to the algebraic Riccati solution in the interval $t = [0,8]$ as $t$ tends to $0$, and it diverges from the infinite horizon solution as $t$ approaches the terminal time.

\section{Future Directions}
There are three immediate future directions for this work. First, the extension to graph systems where each node has multiple states. Second, the extension to graphs which are embedded in metric spaces using vertexon theory \cite{caines_embedded_2022}. So long as the fundamental space the system is embedded in is a Hilbert space, the Q-noise formalism holds. Third, the systems considered in this work have strictly linear noise. This can be expanded into systems of the form of \cite{lu_mathematical_2021}, where the noise intensity depends on the state of the system. 

This paper concerned only centralized control with full observations. One extension would be to apply the Kalman filter to systems with partial observations, and to define a separation principle of the Q-LQG problem and the optimal filtration.



\appendix

\section{Proof of Theorem \ref{thm:mean_convergence}} \label{app:mean_convergence}
    
    Recall 
    \begin{align}
        \bsx^{[N]}_t =& \Phi^{[N]}(t,0) \bsx_0 +\int_{0}^{t} \Phi^{N]}(t,s)   d\bsw^{[N]}_s.
    \end{align}
    Then, 
	\begin{align}
	    &\bE[||\bsx^{[N]}_t - \bsx^{[M]}_t||_2]\\
        &\leq \bE[||\Phi^{[N]}(t,0) \bsx^{[N]}_0 - \Phi^{[M]}(t,0) \bsx^{[M]}_0 ||_2] \nonumber\\ 
        & \quad +\bE[|| \int_{0}^{t} \Phi^{[N]}(t,s)  d\bsw^{[N]}_s - \int_{0}^{t} \Phi^{[M]}(t,s)  d\bsw^{[M]}_s||_2]. \nonumber
	\end{align}

    The Cauchy condition is established in two steps. To address the first expectation, $e^{\bbA^{[N]}t} \bsx^{[M]}$ is added and subtracted inside the norm, which allows assumptions (A0, B1) in conjunction with (A1) to be invoked to bound the term by any $\epsilon'_0 >0$ for sufficiently large $N$ and $M$.

    Next, define $\tbsx_t^{[N]}$ by:
    \begin{equation}
        \tbsx^{[N]}_t = \int_{0}^{t}\Phi^{[N]}(t,s) d\bsw^{[N]}_s.
    \end{equation}
    Then, the second expectation can be evaluated by using the Ito isometry for one-dimensional adapted processes, 
    \begin{equation}
        \bE[||\int_0^T X_t dw_t||^2_2] = \bE[\int_0^T ||X_t ||^2_2 dt].
    \end{equation}
    Then, via the same approach of adding and subtracting $\sum_{k=1}^N e^{\bbA^{[N]}}\bbS^M_k c_{k,r}^M$, the definition of the operator norm gives
    \begin{align}
        &\bE[|| \int_{0}^{t} \Phi^{[N]}(t,s)  d\bsw^{[N]}_s - \int_{0}^{t} \Phi^{[M]}(t,s)  d\bsw^{[M]}_s||_2^2]\\
        & \leq  \int_0^t \bigg(||\Phi^{[N]}(t,s)||_\op^2 \cdot \nonumber \\
        & \quad \quad \quad \quad \sum_{r=1}^M ||(\sum_{i=1}^N \bbS^N_i c_{i,r}^N - \sum_{j=1}^M \bbS^M_j  c_{j,r}^M)||_2^2 \bigg)ds \nonumber\\
        & \quad + \int_0^t \bigg(||\Phi^{[N]}(t,s) - \Phi^{[M]}(t,s)||_\op^2 \cdot  \nonumber \\
        & \quad \quad \quad \quad \quad \sum_{j=1}^M \sum_{r=1}^M ||\bbS^M_j  c_{j,r}^M||_2^2 \bigg)ds\nonumber\\
        & \quad +\int_0^t ||\Phi^{[N]}(t,s)||_\op^2 ds ~ \sum_{r=M+1}^N  \sum_{i=1}^N ||\bbS^N_i||_2^2 |c_{i,r}^N|^2.\nonumber
    \end{align}
    Hence, by the Cauchy property assumptions (A2, A3), the boundedness assumption (B0), and noting that for all $i$, $||\bbS^N_i||^2 = \frac{1}{N^2}$, and by choosing a smaller $\epsilon_2$ if necessary, we see that for any $\epsilon_1' > 0$, and all sufficiently large $N$ and $M$,
    \begin{align}
        &\bE[|| \tbsx^{[N]}_t - \tbsx_t^{[M]}||_2^2]\\
        &<  \bigg(\int_0^t|| \Phi^{[N]}(t,s)||_\op^2 ds ~\epsilon_2 + \epsilon_1  ~\frac{C}{N^2} \nonumber\\
        & \quad\quad\quad\quad\quad\quad\quad + \int_0^t||\Phi^{[N]}(t,s)||_\op^2ds ~ \frac{C}{N^2}\bigg) < \epsilon_1'. \nonumber
    \end{align}
    By the Cauchy-Schwartz inequality applied to the inner product $\langle X, Y \rangle = \bE[XY]$,
    \begin{align}
        \bE[||\tbsx^{[N]} - \tbsx^{[M]}||_2] &=  \sqrt{\bE[1 \cdot ||\tbsx^{[N]} - \tbsx^{[M]}||_2]} \\
        &\leq \sqrt{\bE[1^2]} \sqrt{\bE[||\tbsx^{[N]} - \tbsx^{[M]}||^2_2]} \nonumber\\
        &< \sqrt{\epsilon_1'} \nonumber
    \end{align}
    Hence for any $\epsilon >0$, there exists $N$ and $M$ sufficiently large such that $\sqrt{\epsilon_1'} < \frac{\epsilon}{2}$ and $\epsilon_0 < \frac{\epsilon}{2}$ giving
    $\bE[||\bsx_t^{[N]} - \bsx^{[M]}_t||_2] < \epsilon $. Then, by completeness of $\mathcal{L}^2[0,1]$, there exists $\bsx_t^\infty$ such that $\bE[||\bsx_t^{[N]} - \bsx^\infty_t||_2] \to 0$ as $N$ goes to infinity, yielding the desired result. $\hfill \square$


\section{Proof of Lemma \ref{lemma:uniform_bounded_lemma}} \label{app:uniform_bounded_lemma}
    This lemma is presented in \cite{dunyak_stochastic_2024} and is included here for completeness. 
    Apply the operator norm to both sides of equation (\ref{eqn:uniform_bounded_lemma}). As $\bbR$ (and hence $\bbR^{-1}$) is strictly positive, $-S^N_t \bbB^N \bbR^{-1} \bbB^N S^N_t$ is negative, and by assumption $\{A^N,B^N,Q^N,M^N,M_T^N,R^N\}$ converges to $\{A,B,Q,M,M_T,R\}$ in the operator norm sense, 
    \begin{align}
            ||&{S}^N_t ||_\op \leq  ||\bbM^N_T||_\op + \int_t^T (||\bbA^N S^N_t + S^N_t \bbA^N \nonumber \\
             & \quad \quad \quad  - S^N_t B^N R^{-1} B^N S^N_t + M^N ||_\op) dt \\
            & \leq ||\bbM^N_T||_\op + \int_t^T (4 ||\bbA||_\op ||S^N_t||_\op + 2||\bbM||_\op) dt.  
    \end{align}
    Then, by Gronwall's inequality, $||S^N_t||_\op$ satisfies
    \begin{align}
        ||S^N_t||_\op \leq&  (2||\bbM_T||_\op + 2(T-t) ||\bbM||_\op  ) \nonumber \\
        &~~\cdot \exp(4(T-t) ||\bbA||_\op). 
    \end{align}
    By assumption, $\bbA, \bbM, \bbM_T, \bbR$ are bounded, and hence there exists $0 < c_N < \infty$ such that
    \begin{equation}
        ||S^N_t||_\op \leq 2||\bbM_T||_\op +(T-t)c_N, ~ \forall t \in [0,T], ~ N > N_0. 
    \end{equation}
    $\hfill \square$

\section{Proof of Theorem \ref{thm:LQG_convergence}} \label{app:LQG_convergence}
    This analysis is also presented in \cite{dunyak_stochastic_2024}, but is included here for completeness. Recall that the linear quadratic Q-noise problem is solved by
    \begin{align}
        \bsu_t &= -\bbR^{-1} \bbB^* S_t \bsx_t, ~~ t \in [0,T],\\
        \bsu_t^{[N]} &= -{(\bbR^{[N]})}^{-1} \bbB^{{[N]*}} S^{[N]}_t \bsx^{[N]}_t, ~~t \in [0,T],
    \end{align} which places the system in linear feedback form. Then, by Theorem \ref{thm:mean_convergence}, $\bsx^{[N]}_t$ converges to $\bsx_t$ if
        \begin{equation}
            \lim_{N \to \infty} S^{[N]}_t = S_t, ~~0 \leq t \leq T
        \end{equation}
        in the operator norm sense. Let
        \begin{equation}
            \Delta^N_t := S_t - S^{[N]}_t, ~~ 0\leq t \leq T.
        \end{equation}
        Define the evolution of $\Delta_t^N$ in terms of the evolution of $S_t$ and $S^{[N]}_t$,
        \begin{align}
            \dot{\Delta}^N_t =& \dot{S}_t - \dot{S}_t^{[N]}\\
            =& (\bbA^{[N]} S_t^{[N]} - \bbA S_t) + (S_t^{[N]} \bbA^{[N]*} - S_t \bbA^*) \nonumber\\ 
            &- (S_t^{[N]} \bbB^{[N]} {\bbR^{[N]}}^{-1} \bbB^{[N]*} S_t^{[N]} - S_t\bbB \bbR^{-1} \bbB^* S_t) \nonumber\\
            &  + (\bbM^{[N]} - \bbM), \nonumber \\
            \Delta_T^N =& \bbM_T - \bbM^{[N]}_T,
        \end{align}
        and hence
        \begin{align}
            \Delta^N_t =& \bbM_T - \bbM^{[N]}_T \\
            &+ \int_t^T(\bbA^{[N]} S_t^{[N]} - \bbA S_t) + (S_t^{[N]} \bbA^{[N]*} - S_t \bbA^*) \nonumber\\ 
            &- (S_t^{[N]} \bbB^{[N]} {\bbR^{[N]}}^{-1} \bbB^{[N]*} S_t^{[N]} - S_t\bbB \bbR^{-1} \bbB^* S_t) \nonumber\\
            &+ (\bbM^{[N]} - \bbM)dt. \nonumber
        \end{align}
        Focusing on the first term $ \bbA^{[N]} S_t^{[N]}- \bbA S_t$:
        \begin{align}
            \bbA^{[N]} S_t^{[N]}- \bbA S_t &=\bbA^{[N]} S_t^{[N]}- \bbA S_t+ \bbA^{[N]} S_t - \bbA^{[N]} S_t \nonumber\\
            &= (\bbA^{[N]} - \bbA) S_t - \bbA^{[N]} \Delta^N_t. \label{eqn:lqg_conv1} 
        \end{align}
        Similarly,
        \begin{align}
            S_t^{[N]} \bbA^{[N]*} - S_t \bbA^*   = S_t (\bbA^{[N]*} - \bbA^*) - \Delta^N_t \bbA^{[N]*} \label{eqn:lqg_conv2}
        \end{align}
        For the quadratic term, let 
        \begin{align}
            &H = \bbB \bbR^{-\frac{1}{2}},\\
            &H^{[N]} = \bbB^{[N]} {\bbR^{[N]}}^{-\frac{1}{2}},\\
            &\quad \text{where } \bbR^{-1} = \bbR^{-\frac{1}{2}} \bbR^{-\frac{1}{2}*}. \nonumber
        \end{align}
        Then,
        \begin{align}
            &S_t H H^* S_t = S_t\bbB \bbR^{-1} \bbB^* S_t,\\
            &S_t^{[N]}H^{[N]} H^{[N]*} S_t^{[N]} = S_t^{[N]} \bbB^{[N]} {\bbR^{[N]}}^{-1} \bbB^{[N]*} S_t^{[N]}, \\
            &S_t^{[N]}H^{[N]} H^{[N]*} S_t^{[N]} - S_t H H^* S_t \\
            &\quad = S_t^{[N]}H^{[N]} H^{[N]*} S_t^{[N]} - S_t H H^* S_t\nonumber\\
            & \quad \quad + S_t HH^{[N]*} S_t^{[N]} - S_t HH^{[N]*} S_t^{[N]}\nonumber\\
            & \quad \quad+ S_t H (H^{[N]*}S_t^{[N]} - H^* S_t)\nonumber
        \end{align}
        Employing the identities of the form used in (\ref{eqn:lqg_conv1}) and (\ref{eqn:lqg_conv2}) with $\bbA = H$ and $\bbA^{[N]} = H^{[N]}$, this is equal to
        \begin{align}
            &S_t^{[N]}H^{[N]} H^{[N]*} S_t^{[N]} - S_t H H^* S_t \nonumber \\
            & \quad = \left(S_t (H^{[N]} - H) - \Delta^N_t H^{[N]}\right) H^{[N]*} S_t^{[N]} \nonumber\\
            & \quad \quad + S_t H \left((H^* - H^{[N]*}) S_t - H^{[N]*}\Delta^N_t\right). \label{eqn:lqg_conv3}
        \end{align}
        Using equations (\ref{eqn:lqg_conv1}), (\ref{eqn:lqg_conv2}), (\ref{eqn:lqg_conv3}), define the operator valued functions $P^N,Y^N$,
        \begin{align}
            P^N(t) :=& (\bbM^{[N]} - \bbM) + (\bbA^{[N]} - \bbA) S_t \\
            &+ S_t (\bbA^{[N]*} - \bbA^*) \nonumber\\
            &+S_t (H^{[N]} - H) H^{[N]*} S_t^{[N]} \nonumber\\
            &+ S_t H ( H^* -  H^{[N]*}) S_t, \nonumber
        \end{align}
        \begin{align}            
            Y^N(t,\Delta^N_t) :=& -(\bbA^{[N]}\Delta^N_t + \Delta^N_t \bbA^{[N]*}\\
            & + \Delta^N_t H^{[N]} H^{[N]*} S_t^{[N]} +  S_t H H^{[N]*} \Delta^N_t), \nonumber
        \end{align}
        in terms of which (94) yields
        \begin{equation}
            \Delta_t^N = (\bbM_T - \bbM^{[N]}_T) + \int_t^T P^N(t) + Y^N(t,\Delta_t^N)dt.
        \end{equation}
        Then,
        \begin{align}
            ||\Delta_t^N||_\op  \leq& \int_t^T ||P^N(s) + Y^N(s,\Delta^N_s)||_\op ds \\
            & + ||\bbM_T - \bbM_T^{[N]}||_\op \nonumber\\
            \leq& \int_t^T ||P^N(s)||_\op ds + \int_t^T ||Y^N(s,\Delta^N_s)||_\op ds \nonumber\\
            &+ ||\bbM_T - \bbM_T^{[N]}||_\op. 
        \end{align}
        Further, define the time process $Z^N(t): [0,T] \to \bbR,$
        \begin{align}
            ||Y^N(t,\Delta^N_t)||_\op \leq& (2||\bbA^{[N]}||_\op + ||H^{[N]} H^{[N]*} S_t^{[N]}||_\op \nonumber\\
            &+ ||S_t H H^{[N]*} ||_\op)||\Delta_t^N||_\op \nonumber\\
            =: & Z^N(t) ||\Delta_t^N||_\op. \label{eq:app_deltat}
        \end{align}
        Then, by applying Gronwall's inequality to $||\Delta_t^N||_\op$ in (\ref{eq:app_deltat}),
        \begin{align} \label{eq:delta_gronwall}
            ||\Delta_t^N||_\op \leq& (||\bbM_T - \bbM_T^{[N]}||_\op  \\
            &+ \int_t^T ||P^N(s)||_\op ds) \exp( \int_t^T ||Z^N(s)||_\op ds). \nonumber
        \end{align}
        By assumption $\bbA^{[N]} \to \bbA,~H^{[N]} \to H,~ \bbM^{[N]} \to \bbM, ~ \bbM^{[N]}_T \to \bbM_T$ in operator norm, and as $S_t, S^{[N]}_t$ are uniformly bounded operators for all $N, t$,
        \begin{align} 
            ||&P^N(t)||_\op\\
            & = ||(\bbM^{[N]} - \bbM) + (\bbA^{[N]} - \bbA) S_t \nonumber \\
            &+ S_t (\bbA^{[N]*} - \bbA^*) \nonumber \\
            &+S_t (H^{[N]} - H) H^{[N]*} S_t^{[N]} \nonumber \\
            &+ S_t H ( H^* -  H^{[N]*}) S_t||_\op\nonumber \\
            \leq& ||(\bbM^{[N]} - \bbM)||_\op + ||(\bbA^{[N]} - \bbA)||_\op ||S_t||_\op \label{eq:P_convergence} \\
            &+||S_t||_\op ||(\bbA^{[N]*} - \bbA^*)||_\op \nonumber\\
            &+||S_t||_\op ||(H^{[N]} - H)||_\op ||H^{[N]*}||_\op ||S_t^{[N]}||_\op \nonumber \\
            &+ ||S_t||_\op ||H||_\op ||( H^* -  H^{[N]*})||_\op ||S_t||_\op. \nonumber 
        \end{align}
        Hence, by the convergence of $\bbM^{[N]}_T \to \bbM_T$ and equation (\ref{eq:P_convergence}), as $N \to \infty$,
        \begin{align}
            &(||\bbM_T - \bbM_T^{[N]}||_\op + \int_t^T ||P^N(s)||_\op ds) \to 0\\
            &\text{and } ||Z(t)||_\op < \infty \implies \exp(\int_t^T ||Z(s)||_\op ds) < \infty.
        \end{align}
        Hence, by (\ref{eq:delta_gronwall})
        \begin{align}
            ||\Delta^N_t||_\op \to 0.
        \end{align}
        As $||\Delta_t^N||_\op$ converges to zero, $S^{[N]}_t$ converges to $S_t$ in the operator norm sense as $N$ increases to infinity. Then the finite dimensional network operator $(\bbA^{[N]} - \bbB^{[N]} \bbK^{[N]}_t)$ converges to the operator on the graphon system $(\bbA - \bbB \bbK_t)$ and Theorem \ref{thm:mean_convergence} can be applied, yielding
        \begin{equation}
            \bE[||\bsx_t^{[N]} - \bsx_t||_2] < \epsilon,
        \end{equation}
        as required.
        $\hfill \square$
\bibliographystyle{plain}
\bibliography{Q_graphon_auto.bib}

\end{document}